\documentclass[12pt]{amsart}
\usepackage{amscd,amssymb,tikz,pgflibraryarrows}
\newcommand\field[1]{\mathbb{#1}}
\newcommand\NN{\field{N}}
\newcommand\QQ{\field{Q}}
\newcommand\TT{\field{T}}
\newcommand\ZZ{\field{Z}}
\newcommand\Kk{\mathcal K}

\newcommand\Mm{\mathcal M}

\renewcommand\ker{\operatorname{ker}}
\newcommand\coker{\operatorname{coker}}

\newcommand\id{\operatorname{id}}

\newcommand\Aut{\operatorname{Aut}}
\newcommand\lsp{\operatorname{span}}
\newcommand\clsp{\operatorname{\overline{span\!}\,\,}}
\newcommand\MCE{\operatorname{MCE}}

\newcommand\im{\operatorname{im}}
\newcommand\col{\operatorname{col}}

\newcommand{\cpg}[3]{{#1 \times_{#2} \ZZ^{#3}}}
\newcommand{\Z}[1]{0_{#1}}
\newcommand{\indicator}{1}
\newcommand{\tgrphlim}{{
    \renewcommand{\leftarrow}{\leftharpoondown}
    \varprojlim
}}

\theoremstyle{plain}
\newtheorem{theorem}{Theorem}[section]
\newtheorem*{theorem*}{Theorem}
\newtheorem*{prop*}{Proposition}
\newtheorem{cor}[theorem]{Corollary}
\newtheorem{lemma}[theorem]{Lemma}
\newtheorem{prop}[theorem]{Proposition}

\theoremstyle{remark}
\newtheorem{rmk}[theorem]{Remark}

\newtheorem{example}[theorem]{Example}
\newtheorem{examples}[theorem]{Examples}
\theoremstyle{definition}
\newtheorem{dfn}[theorem]{Definition}

\newtheorem{notation}[theorem]{Notation}

\numberwithin{equation}{section}

\newcounter{myenumi}
\newenvironment{Enumerate}
{\begin{list} {\textnormal{(\arabic{myenumi})}}
{\usecounter{myenumi}\setlength{\leftmargin}{1.5em}\setlength{\labelwidth}{2.5em}}}
{\end{list}}

\title{Crossed products of $k$-graph $C^*$-algebras by $\ZZ^l$.}
\author{Cynthia Farthing}
\address{Cynthia Farthing \\
University of Nebraska-Lincoln \\
Department of Mathematics \\
203 Avery Hall \\
Lincoln\\ NE 68588-0130 \\ USA}
\email{cfarthing2@unl.edu}
\author{David Pask}
\address{David Pask \\
School of Mathematics and Applied Statistics \\
Austin Keane Building (15) \\
University of Wollongong \\
NSW 2522 \\ AUSTRALIA}
\email{dpask@uow.edu.au}
\author{Aidan Sims}
\address{Aidan Sims \\
School of Mathematics and Applied Statistics \\
Austin Keane Building (15) \\
University of Wollongong \\
NSW 2522 \\ AUSTRALIA}
\email{asims@uow.edu.au}
\date{June 25, 2007}
\thanks{This research was supported by the Australian Research
Council.}

\begin{document}

\begin{abstract}
An action of $\ZZ^l$ by automorphisms of a $k$-graph induces an
action of $\ZZ^l$ by automorphisms of the corresponding
$k$-graph $C^*$-algebra. We show how to construct a
$(k+l$)-graph whose $C^*$-algebra coincides with the crossed
product of the original $k$-graph algebra by $\ZZ^l$. We then
investigate the structure of the crossed-product $C^*$-algebra.
\end{abstract}

\maketitle

\section{Introduction}

In recent years, much attention has been paid to graph algebras
and their higher-rank analogues as models for classifiable
$C^*$-algebras (see \cite{CBMSbk} for an overview of the
subject). What makes these models so attractive is the ability
to trade information back and forth between the underlying
combinatorial object and the associated $C^*$-algebra. This
program is quite advanced for graph algebras (see, for example,
\cite{HS2, Jeong, KPR, Sz2}). However, higher-rank graphs are a
more recent development and have more complicated combinatorial
properties than ordinary graphs. Consequently, many important
structural questions regarding higher-rank graph algebras
remain unanswered. In particular, while the general theory of
$k$-graph $C^*$-algebras is quickly catching up with that of
graph $C^*$-algebras (see for example \cite{Evans, Far, KP,
RSY2, RobSi}), there remains a dearth of tractable examples in
the higher-rank setting.

A construction which \emph{has} been successfully generalised
from the setting of graphs to that of higher-rank graphs is the
skew-product construction (see \cite{KP, PQR2}). This
construction allows us to realise certain crossed products of
$k$-graph algebras by coactions of groups $G$ as $k$-graph
algebras in their own right. Specifically, suppose that the
coaction $\delta$ arises from a functor $c$ from the $k$-graph
$\Lambda$ to the group $G$ (that is, a function which takes
composition in $\Lambda$ to multiplication in $G$). Then we may
form a skew-product $k$-graph $\Lambda \times_c G$, and the
coaction crossed product $C^*(\Lambda) \times_\delta G$ is
canonically isomorphic to $C^*(\Lambda \times_c G)$.

Results of~\cite{KP1, KP} show that we can also realise certain
crossed products of $k$-graph $C^*$-algebras by \emph{actions}
of groups as $k$-graph $C^*$-algebras. This is achieved using
the quotient-graph construction. Specifically, an action
$\alpha$ of a group $G$ by automorphisms of a $k$-graph
$\Lambda$ induces an action $\tilde{\alpha}$ of $G$ by
automorphisms of $C^*(\Lambda)$. If $\alpha$ is free in the
sense that no nontrivial group element fixes any vertex, then
one may form the quotient $k$-graph $\Lambda/G$ whose morphisms
are the orbits of morphisms of $\Lambda$ under the action of
$G$. In this situation, the $C^*$-algebras $C^*(\Lambda)
\times_{\tilde\alpha} G$ and $C^*(\Lambda/G)$ are Morita
equivalent \cite[Theorem~5.7]{KP}. If the action $\alpha$ is
not free, however, \cite[page~176]{PQR2} shows that $\Lambda/G$
may not be a $k$-graph because composition may not be
well-defined, so the approach of \cite{KP1, KP} is not
applicable. Moreover, notice that whereas there is an
isomorphism between a skew-product $k$-graph $C^*$-algebra and
the associated coaction crossed product $C^*$-algebra, the
corresponding result for actions using the quotient $k$-graph
construction yields only a Morita equivalence.

In this article we describe a class of $(k+l)$-graphs whose
$C^*$-algebras are isomorphic to crossed products of $k$-graph
algebras by $\ZZ^l$. Given an action $\alpha$ of $\ZZ^l$ on a
$k$-graph $\Lambda$, we construct a $(k+l)$-graph
$\cpg{\Lambda}{\alpha}{l}$ which, as the notation suggests, can
profitably be thought of as the crossed product of $\Lambda$ by
$\alpha$. We show that $\alpha$ induces an action
$\tilde{\alpha}$ of $\ZZ^l$ on $C^*(\Lambda)$ and that the
higher-rank graph $C^*$-algebra $C^*(\cpg{\Lambda}{\alpha}{l})$
coincides with the crossed-product $C^*$-algebra $C^*(\Lambda)
\times_{\tilde{\alpha}} \ZZ^l$.

It is noteworthy that our results do not require that $\alpha$
should be free, and we obtain an isomorphism rather than a
Morita equivalence. Moreover, there is good evidence to suggest
that our notion of a crossed-product $k$-graph is a reasonable
one. For example, we combine our construction with the
skew-product construction of \cite{KP, PQR2} to obtain a
realisation of Takai duality at the level of higher-rank graphs
(see Section~\ref{sec:takai}).

The identification of $C^*(\cpg{\Lambda}{\alpha}{l})$ with
$C^*(\Lambda) \times_{\tilde\alpha} \ZZ^l$ allows us to study
the crossed-product $C^*$-algebra using the theory of graph
$C^*$-algebras, and in particular to formulate necessary and
sufficient conditions for simplicity of the crossed product. It
also allows us to study the $(k+l)$-graph $C^*$-algebra using
the theory of crossed-product $C^*$-algebras; for example, when
$l=1$, our construction lends itself to analysis of the
$K$-theory of $C^*(\Lambda \times_\alpha \ZZ)$ via the
Pimsner-Voiculescu exact sequence.

\medskip

The paper is organised as follows. In
Section~\ref{sec:prelims}, we introduce the notation and
conventions we will use throughout. We open
Section~\ref{sec:the construction} by describing our
construction of a $(k+l)$-graph $\cpg{\Lambda}{\alpha}{l}$ from
an action $\alpha$ of $\ZZ^l$ on a $k$-graph $\Lambda$. We
establish that the action $\alpha$ of $\ZZ^l$ on $\Lambda$
induces an action $\tilde{\alpha}$ of $\ZZ^l$ on
$C^*(\Lambda)$, and then prove our first main result,
Theorem~\ref{thm:crossed product}: the $(k+l)$-graph algebra
$C^*(\cpg{\Lambda}{\alpha}{l})$ is canonically isomorphic to
the crossed product $C^*$-algebra $C^*(\Lambda)
\times_{\tilde{\alpha}} \ZZ^l$.

In Section~\ref{sec:simplicity}, we use the recent results of
\cite{RobSi} which characterise simplicity of higher-rank graph
algebras to decide when $C^*(\Lambda \times_{\alpha} \ZZ^l)$ is
simple in terms of properties of $\Lambda$ and the action
$\alpha$. In Section~\ref{sec:C*-simplicity}, we recast the
results of Section~\ref{sec:simplicity} in terms of features of
the $C^*$-algebra $C^*(\Lambda)$ and of properties of the
induced action $\tilde{\alpha}$ of $\ZZ^l$ on $C^*(\Lambda)$.

We conclude in section~\ref{sec:K-theory} with an application
of our results to the calculation of $K$-theory for certain
examples. Specifically, we consider a $1$-graph $E$ endowed
with an action $\alpha$ of $\ZZ$ so that $\tilde{\alpha}$ is an
action of $\ZZ$ on $C^*(E)$. If either of the $K$-groups of
$C^*(E)$ is trivial, we may use the Pimsner-Voiculescu exact
sequence in $K$-theory to calculate the $K$-groups of
$C^*(\cpg{E}{\alpha}{})$.

\medskip

\textbf{Acknowledgements.} The authors wish to thank Iain Raeburn and
Astrid an Huef for a number of helpful discussions.

\section{Preliminaries}\label{sec:prelims}
We first recall the notation and conventions used for
$k$-graphs. For more details see \cite{CBMSbk, RSY2}.

\subsection{The semigroup $\NN^k$.}
We write $\NN$ for the semigroup $\{0,1,2,\dots\}$ under
addition. We regard $\NN^k$ as a semigroup under addition with
identity element denoted $\Z{k}$. When convenient, we regard
the semigroup $\NN^k$ as (the morphisms of) a category with a
single object (the addition operator on $\NN^k$ is viewed as a
composition map).

We denote the canonical generators of $\NN^k$ by $e_1, \dots,
e_k$, and for $n \in \NN^k$ and $1 \le i \le k$, we write $n_i$
for the $i^{\rm th}$ coordinate of $n$, so that $n =
\sum^k_{i=1} n_i e_i$. Fix $m,n \in \NN^k$. We write $m \le n$
if $m_i \le n_i$ for all $i$. We denote by $m \vee n$ the
coordinate-wise maximum of $m$ and $n$, and $m \wedge n$ the
coordinatewise minimum. In particular, we have $m \wedge n \le
m,n \le m \vee n$.

We shall often and without comment identify $\NN^{k+l}$ with $\NN^k
\times \NN^l$. In particular, we write $(p,m) \in \NN^{k+l}$ to
indicate that $(p,m)$ is the element of $\NN^{k+l}$ whose first $k$
coordinates are those of $p \in \NN^k$ and whose last $l$ coordinates
are those of $m \in \NN^l$.

\subsection{Higher-rank graphs.}
Recall from \cite[Definition~1.1]{KP} that a $k$-graph is a
countable category $\Lambda$ together with a functor $d :
\Lambda \to \NN^k$ which satisfies the factorisation property:
if $\lambda \in \Lambda$ with $d(\lambda) = m + n$, then there
exist unique elements $\mu \in d^{-1}(m)$ and $\nu \in
d^{-1}(n)$ such that $\lambda = \mu\nu$. We call $d$ the
\emph{degree functor} and regard it as a higher-rank analogue
of length. An argument involving the factorisation property
shows that $v \mapsto \id_v$ is a bijection between the objects
of $\Lambda$ and the morphisms of degree $\Z{k}$. We use this
to identify the two, and we regard $\Lambda$ as a collection of
morphisms only.

For $n \in \NN^k$, we denote $d^{-1}(n)$ by $\Lambda^n$. We
call the elements of $\Lambda$ \emph{paths} and the elements of
$\Lambda^{\Z{k}}$ \emph{vertices}. If $\lambda \in \Lambda^p$
and $\Z{k} \le m \le n \le p$ then we denote by
$\lambda(\Z{k},m)$, $\lambda(m,n)$ and $\lambda(n,p)$ the
unique elements of $\Lambda^m$, $\Lambda^{n-m}$ and
$\Lambda^{p-n}$ satisfying $\lambda =
\lambda(\Z{k},m)\lambda(m,n)\lambda(n,p)$ (the existence and
uniqueness of these morphisms follows from two applications of
the factorisation property.)

For $\mu,\nu \in \Lambda$ we call $\lambda$ a \emph{common
extension} of $\mu$ and $\nu$ if $\lambda = \mu\alpha =
\nu\beta$ for some $\alpha, \beta \in \Lambda$. If $\lambda$ is
a common extension of $\mu$ and $\nu$ then we must have
$d(\lambda) \ge d(\mu) \vee d(\nu)$. We call $\lambda$ a
\emph{minimal common extension} of $\mu$ and $\nu$ if it is a
common extension satisfying $d(\lambda) = d(\mu) \vee d(\nu)$.
We write $\MCE(\mu,\nu)$ for the set of all minimal common
extensions of $\mu$ and $\nu$. We say that $\Lambda$ is
\emph{finitely aligned} if $\MCE(\mu,\nu)$ is finite (possibly
empty) for all $\mu,\nu \in \Lambda$. We write $\Lambda^{\rm
min}(\mu,\nu)$ for the set $\{(\alpha,\beta) : \mu\alpha =
\nu\beta \in \MCE(\mu,\nu)\}$.

If $S$ is a subset of $\Lambda$ and $v \in \Lambda^{\Z{k}}$, we write
$vS$ for the set $S \cap r^{-1}(v)$ and we write $Sv$ for $S \cap
s^{-1}(v)$. We say that $\Lambda$ is \emph{row-finite} if
$|v\Lambda^n| < \infty$ for all $v \in \Lambda^{\Z{k}}$ and $n \in
\NN^k$, and we say that $\Lambda$ has \emph{no sources} if
$v\Lambda^n \not= \emptyset$ for all $v \in \Lambda^{\Z{k}}$ and $n
\in \NN^k$.

Given a vertex $v \in \Lambda^{\Z{k}}$ and a subset $F$ of
$v\Lambda$, we say $F$ is \emph{exhaustive} if for every $\mu
\in v\Lambda$ there exists $\nu \in F$ such that $\MCE(\mu,\nu)
\not= \emptyset$. We say an exhaustive set $F \subset v\Lambda$
is \emph{finite exhaustive} if $|F| < \infty$. If $\Lambda$ has
no sources, then $v\Lambda^n$ is exhaustive for all $v \in
\Lambda^{\Z{k}}$ and $n \in \NN^k$, so if $\Lambda$ is
row-finite and has no sources, then $v\Lambda^n$ is always
finite exhaustive.

\subsection{The universal $C^*$-algebra $C^*(\Lambda)$.}
As in \cite{RSY2}, given a finitely aligned $k$-graph
$\Lambda$, a Cuntz-Krieger $\Lambda$-family is a set
$\{t_\lambda : \lambda \in \Lambda\}$ of partial isometries
satisfying
\begin{itemize}
\item[(TCK1)]
$\{t_v : v \in \Lambda^{\Z{k}}\}$ is a set of mutually orthogonal
projections.
\item[(TCK2)]
$t_\mu t_\nu = t_{\mu\nu}$ whenever $r(\nu) = s(\mu)$.
\item[(TCK3)]
$t^*_\mu t_\nu = \sum_{(\xi,\eta) \in \Lambda^{\rm min}(\mu,\nu)}
t_\xi t^*_\eta$.
\item[(CK)] $\prod_{\lambda \in F} (t_v - t_\lambda
    t^*_\lambda) = 0$ for each $v \in \Lambda^{\Z{k}}$ and
    each finite exhaustive $F \subset v\Lambda$.
\end{itemize}
By $t_\Lambda$, we mean the Cuntz-Krieger $\Lambda$-family
$\{t_\lambda : \lambda \in \Lambda\}$ as a whole.

General results of Blackadar (see \cite{Bla}) imply that there
is a $C^*$-algebra $C^*(\Lambda)$ (unique up to canonical
isomorphism) generated by a Cuntz-Krieger $\Lambda$-family
$s_\Lambda$ which is universal in the sense that for any other
Cuntz-Krieger $\Lambda$-family $t_\Lambda$ there is a
homomorphism $\pi_t : C^*(\Lambda) \to C^*(t_\Lambda)$
satisfying $\pi_t(s_\lambda) = t_\lambda$ for all $\lambda \in
\Lambda$.

For $z \in \TT^k$ and $n \in \ZZ^k$, we employ multi-index
notation and write $z^n$ for the product $\prod_{i=1}^k
z_i^{n_i} \in \TT$. Using the universal property of
$C^*(\Lambda)$ one can check that there is a strongly
continuous action $\gamma$ of $\TT^k$ on $C^*(\Lambda)$
satisfying $\gamma_z(s_\lambda) = z^{d(\lambda)}s_\lambda$ for
all $\lambda \in \Lambda$ and $z \in \TT^k$. This action
$\gamma$ is called the \emph{gauge action}.

\subsection{Graph morphisms and infinite
paths}%\label{sec:alpha^infty}

Given $k$-graphs $\Lambda$ and $\Gamma$, a \emph{$k$-graph
morphism} $\phi : \Lambda \to \Gamma$ is a functor from
$\Lambda$ to $\Gamma$ which respects the degree maps. A
bijective $k$-graph morphism is simply called an isomorphism.
An isomorphism $\phi : \Lambda \to \Lambda$ is called an
automorphism of $\Lambda$.

To discuss infinite paths in $k$-graphs, we must first
introduce the $k$-graph $\Omega_k$. For $k \ge 1$, $\Omega_k$
is the $k$-graph $\{(m,n) \in \NN^k \times \NN^k : m \le n\}$
with structure maps $r(m,n) := (m,m)$, $s(m,n) := (n,n)$,
$(m,n)(n,p) := (m,p)$ and $d(m,n) := n-m$. We typically denote
the element $(m,m)$ of $\Omega_k^{\Z{k}}$ by $m$.

An \emph{infinite path} in a $k$-graph $\Lambda$ is a graph
morphism $x : \Omega_k \to \Lambda$. We denote the collection
of all infinite paths in $\Lambda$ by $\Lambda^\infty$. For any
cofinal sequence $(m_i)^\infty_{i=1} \subset \NN^k$ such that
$m_0 = \Z{k}$, a pair of infinite paths $x,y$ are equal if and
only if $x(m_i, m_{i+1}) = y(m_i, m_{i+1})$ for all $i$. Hence,
given a cofinal sequence $(m_i)^\infty_{i=1} \subset \NN^k$, we
may view an infinite path $x$ as the infinite composition of
finite paths $x = x(0_k,m_1) x(m_1, m_2) x(m_2,m_3) \cdots$. In
keeping with this, we regard $x(\Z{k})$ as the range of $x$ and
denote it $r(x)$.

Fix $x \in \Lambda^\infty$. For each $\lambda \in \Lambda r(x)$
there is a unique infinite path $\lambda x$ satisfying
$(\lambda x)(\Z{k}, d(\lambda)) = \lambda$ and $(\lambda
x)(d(\lambda), d(\lambda) + m) = x(\Z{k}, m)$ for all $m \in
\NN^k$. We write $\lambda \Lambda^\infty$ for the set
$\{\lambda x : x \in \Lambda^\infty, r(x) = s(\lambda)\}
\subset \Lambda^\infty$, and call this the \emph{cylinder set
associated to $\lambda$}. For each $m \in \NN^k$, there is a
unique infinite path $\sigma^m(x)$ such that
$\sigma^m(x)(\Z{k},n) = x(m, m+n)$ for all $n \in \NN^k$. Note
that $\sigma^{d(\lambda)}(\lambda x) = x = x(0,m)\sigma^m(x)$
for each $\lambda \in \Lambda r(x)$ and each $m \in \NN^k$.

\subsection{The $k$-graph $\Delta_k$}\label{sec:Delta k}
Related to infinite paths and the $k$-graph $\Omega_k$ is the
two-sided version $\Delta_k$ of $\Omega_k$. For an integer $k
\ge 1$, we write $\Delta_k$ for the $k$-graph $\Delta_k :=
\{(m,n) \in \ZZ^k \times \ZZ^k : m \le n\}$ with $r(m,n) = m$,
$s(m,n) = n$, $d(m,n) = n-m$ and $(m,n)(n,p) := (m,p)$.

Note that $\Omega_k$ is isomorphic to the sub-$k$-graph of
$\Delta_k$ consisting of elements $(m,n)$ such that $m \ge 0$.
As with $\Omega_k$ we usually denote a vertex $(m,m)$ of
$\Delta_k$ by $m$.

\subsection{Skeletons.}\label{sec:skeletons}
One can completely describe a $k$-graph using its skeleton,
which consists of a $k$-coloured graph $E$ and a list of
factorisation rules $F$. We outline this construction here, but
see \cite[pp. 90--91]{CBMSbk} or \cite[Section~2]{RSY1} for
more detail.

By a $k$-coloured graph, we mean a $5$-tuple $(E^0, E^1, r, s,
\col)$ where $(E^0, E^1, r, s)$ is a directed graph, and $\col
: E^1 \to \{1, \dots, k\}$ is the \emph{colour} function. For
convenience, we denote $\col^{-1}(i) \subset E^1$ by $E^1_i$.

To each $k$-graph $\Lambda$ we associate a $k$-coloured graph
$E_\Lambda = (\Lambda^{\Z{k}}, \bigsqcup^k_{i=1} \Lambda^{e_i},
r, s, \col)$ where the range and source maps $r,s$ are
inherited from $\Lambda$, and where $\col(f) = i$ if and only
if $f \in \Lambda^{e_i}$. In the examples in this article, two
colours will suffice: edges of degree $e_1$ will be thought of
as blue (and drawn using solid lines), and edges of degree
$e_2$ will be thought of as red (and drawn using dashed lines).

Let $1 \le i < j \le k$, and suppose $f,g \in E_\Lambda^1$ with
$\col(f) = i$, $\col(g) = j$ and $s(f) = r(g)$. Then $fg \in
\Lambda^{e_i + e_j}$, so the factorisation property in
$\Lambda$ applied with $d(fg) = e_j + e_i$ ensures that there
are unique edges $g', f' \in E_\Lambda^1$ such that $\col(f') =
i$, $\col(g') = j$, $s(g') = r(f')$, and $fg = g'f'$ in
$\Lambda$. The list of factorisation rules $F_\Lambda$
associated to $\Lambda$ is the complete set of equalities $f g
= g' f'$ obtained this way. The skeleton of $\Lambda$ is the
pair $(E_\Lambda, F_\Lambda)$.

Conversely, let $E = (E^0, E^1, r, s, \col)$ be a $k$-coloured
directed graph. Let $F$ be a list of equalities of the form $f
g = g'f'$ where $f,f' \in E^1_i$, $g,g' \in E^1_j$, $i < j$,
$s(f) = r(g)$ and $s(g') = r(f')$. We say that $F$ is
\emph{permissible} it satisfies two conditions. The first
condition is essentially the factorisation property for
bi-coloured paths of length~2:
\begin{itemize}
\item[(1)]
the factorisation rules determine bijections $E^1_i \times_{E^0}
E^1_j \to E^1_j \times_{E^0} E^1_i$ for $i < j$, where $E^1_i
\times_{E^0} E^1_j$ is the fibred product $\{(f,g) \in E^1_i
\times E^1_j : s(f) = r(g)\}$. That is, each bi-coloured path in
$E$ appears in exactly one factorisation rule.
\end{itemize}
To state the second rule, observe that if $(f,g,h) \in E^1_i
\times_{E^0} E^1_j \times_{E^0} E^1_l$ (where $i,j,l \in \{1,
\dots, k\}$ are distinct), then~(1) gives unique edges $f^1,
f^2 \in E^1_i$, $g^1, g^2 \in E^1_j$ and $h^1, h^2 \in E^1_k$
such that
\[
fgh = fh^1g^1 = h^2f^1g^1 = h^2f^2g^2.
\]
Likewise there are unique edges $f_1, f_2 \in E^1_i$, $g_1, g_2 \in
E^1_j$ and $h_1, h_2 \in E^1_k$ such that
\[
fgh = g_1f_1h = g_1h_1f_2 = h_2g_2f_2.
\]
The collection $F$ of factorisation rules is permissible if it
satisfies~(1) and
\begin{itemize}
\item[(2)] The factorisation rules are associative: for
    every $(f,g,h) \in E^1_i \times_{E^0} E^1_j
    \times_{E^0} E^1_l$ such that $i,j,l \in \{1, \dots,
    k\}$ are distinct, the two different ways of reversing
    the colours in the path $fgh$ discussed above agree:
    $f^2 = f_2$, $g^2 = g_2$ and $h^2 = h_2$.
\end{itemize}
The pair $(E,F)$ is then called a skeleton. The results of
\cite{FS1} imply that the map $\Lambda \mapsto (E_\Lambda,
F_\Lambda)$ which sends a $k$-graph $\Lambda$ to its skeleton
is reversible: given a skeleton $(E,F)$, there is a unique
$k$-graph $\Lambda_{(E,F)}$ such that $(E_{\Lambda_{(E,F)}},
F_{\Lambda_{(E,F)}}) = (E,F)$.

Note that if $k = 1$, then there are no factorisation rules to
list, so every ($1$-coloured) directed graph specifies a
$1$-graph. Likewise, if $k = 2$, then~(2) above is trivial
because we cannot have distinct $i,j,l \in \{1, \dots, k\}$, so
every bi-coloured graph together with factorisation rules
satisfying~(1) specifies a $2$-graph.

\section{Crossed products by $\ZZ^l$ as higher-rank graph
algebras}\label{sec:the construction}

In this section we show how an action $\alpha$ of $\ZZ^l$ on a
finitely aligned $k$-graph $\Lambda$ induces an action
$\tilde{\alpha}$ of $\ZZ^l$ on $C^*(\Lambda)$, generalising the
assertion of \cite[page~16]{KP} to the finitely aligned
setting. We then show how to realise the crossed-product
$C^*$-algebra $C^*(\Lambda) \times_{\tilde{\alpha}} \ZZ^l$ as
the universal algebra of a $(k+l)$-graph
$\cpg{\Lambda}{\alpha}{l}$.

\begin{prop}\label{prp:obtain action}
Let $\Lambda$ be a finitely aligned $k$-graph and let
$s_\Lambda$ be the universal Cuntz-Krieger $\Lambda$-family in
$C^*(\Lambda)$.
\begin{enumerate}
\item Let $\phi$ be an automorphism of $\Lambda$. Then
    there is a unique automorphism $\tilde{\phi}$ of
    $C^*(\Lambda)$ satisfying $\tilde{\phi}(s_\lambda) =
    s_{\phi(\lambda)}$ for all $\lambda \in \Lambda$.
\item Let $G$ be a group, and suppose that $g \mapsto
    \alpha_g$ is an action of $G$ on $\Lambda$ by
    automorphisms. Then $g \mapsto \tilde\alpha_g$ is an
    action of $G$ on $C^*(\Lambda)$ by automorphisms.
\end{enumerate}
\end{prop}
\begin{proof}
1) It is easy to check that a graph isomorphism preserves
minimal common extensions and finite exhaustive sets. It
follows from this that there is a Cuntz-Krieger
$\Lambda$-family $t^\phi_\Lambda$ defined by
$t^\phi_\lambda := s_{\phi(\lambda)}$ for $\lambda \in
\Lambda$. The universal property of $C^*(\Lambda)$
therefore furnishes us with a $C^*$-homomorphism
$\tilde{\phi}$ of $C^*(\Lambda)$ satisfying
$\tilde{\phi}(s_\lambda) = t^\phi_\lambda
=s_{\phi(\lambda)}$ for all $\lambda \in \Lambda$. Applying
the same argument to $\phi^{-1} \in \Aut(\Lambda)$ gives
another $C^*$-homomorphism $\widetilde{\phi^{-1}} :
C^*(\Lambda) \to C^*(\Lambda)$, and since
$\widetilde{\phi^{-1}} \circ \tilde{\phi}$ fixes all the
generators of $C^*(\Lambda)$, $\widetilde{\phi^{-1}}$ is an
inverse for $\tilde\phi$ and in particular, $\tilde\phi$ is
an automorphism.

2) Let $1_G$ denote the identity element of $G$. Since
$\alpha$ is an action,
\[
 \tilde{\alpha}_{1_G}(s_\lambda) = s_\lambda,\quad
 \tilde{\alpha}_{g^{-1}}(\tilde{\alpha}_g(s_\lambda)) = s_\lambda
 \quad\text{and}\quad
 \tilde{\alpha}_g(\tilde{\alpha}_h(s_\lambda)) =
 \tilde{\alpha}_{gh}(s_\lambda),
\]
for each generator $s_\lambda$ of $C^*(\Lambda)$ and all $g,h
\in G$. It follows that $\tilde{\alpha}$ is an action of
$\ZZ^l$ on $C^*(\Lambda)$ by automorphisms.
\end{proof}

We now show how to construct a $(k+l)$-graph
$\cpg{\Lambda}{\alpha}{l}$ from an action of $\ZZ^l$ on a
$k$-graph $\Lambda$. We show in Theorem~\ref{thm:crossed
product} that the $C^*$-algebra of this $(k+l)$-graph is
isomorphic to the crossed-product $C^*(\Lambda)
\times_{\tilde{\alpha}} \ZZ^l$.

\begin{prop}\label{prp:Lambda^alpha}
Let $\Lambda$ be a $k$-graph, and suppose that $\alpha$ is an
action of $\ZZ^l$ on $\Lambda$ by automorphisms. Then there is
a unique $(k+l)$-graph $\cpg{\Lambda}{\alpha}{l}$ such that
\begin{Enumerate}
\item $(\cpg{\Lambda}{\alpha}{l})^{(p,m)} = \Lambda^p \times \{m\}$
for all $(p,m) \in \NN^{k+l}$;
\item $r(\lambda,m) = (r(\lambda),\Z{l})$ and $s(\lambda,m) =
(\alpha_{-m}(s(\lambda)), \Z{l})$ for all $\lambda \in \Lambda$ and
$m \in \NN^l$; and
\item $(\mu,m)(\nu,n) = (\mu\alpha_m(\nu), m+n)$ whenever $s(\mu,m)
    =
r(\nu,n)$.
\end{Enumerate}
Moreover, $\cpg{\Lambda}{\alpha}{l}$ is finitely aligned if and
only if $\Lambda$ is finitely aligned,
$\cpg{\Lambda}{\alpha}{l}$ is row-finite if and only if
$\Lambda$ is row-finite, and $\cpg{\Lambda}{\alpha}{l}$ has no
sources if and only if $\Lambda$ has no sources.
\end{prop}

\begin{rmk}
Although the notation may suggest otherwise, the $k$-graph
$\cpg{\Lambda}{\alpha}{l}$ is equal as a set to $\Lambda \times
\NN^l$ rather than to $\Lambda \times \ZZ^l$.
\end{rmk}

\begin{proof}[Proof of Proposition~\ref{prp:Lambda^alpha}]
The details of this proof are quite messy, but the idea is
straightforward; we present only the idea here.

The discussion in Section~\ref{sec:skeletons} shows that it
suffices to produce the skeleton of $\cpg{\Lambda}{\alpha}{l}$.
That is, a $(k+l)$-coloured graph $E$ and an allowable list $F$
of factorisation rules. To obtain $E$, we begin with a copy
$E_\Lambda \times \{\Z{l}\}$ of the $k$-coloured graph
associated to $\Lambda$ and augment it as follows: for each $1
\le i \le l$ and each $v \in \Lambda^{\Z{k}}$, we add an edge
$(v, e_i)$ to $E_{k+i}^1$ directed from the vertex
$(\alpha_{e_i}^{-1}(v),\Z{l})$ to the vertex $(v,\Z{l})$. The
factorisation rules $F$ are specified as follows.
\begin{itemize}
\item If $f_1 f_2 = f'_2 f'_1$ in $\Lambda$, then $(f_1,
    \Z{l})(f_2, \Z{l}) = (f'_2, \Z{l}) (f'_1, \Z{l})$
    belongs to $F$; that is, the factorisation rules
    amongst edges from $E_\Lambda$ are unchanged.
\item For $f \in \Lambda^{e_j}$ and $1 \le i \le l$,
    $(f,\Z{l}) (s(f), e_i) = (r(f), e_i)
    (\alpha_{e_i}^{-1}(f),\Z{l})$.
\item For $i \not= j$ in $\{1, \dots, l\}$ and $v \in
    \Lambda^{\Z{k}}$, $(v, e_i)(\alpha_{e_i}^{-1}(v), e_j)
    = (v, e_j) (\alpha_{e_j}^{-1}(v), e_i)$.
\end{itemize}
We must check that this collection $F$ satisfies conditions
(1)~and~(2) of Section~\ref{sec:skeletons}. That $\alpha$ is an
action ensures that the $\alpha_{e_i}$ commute, so
Condition~(1) is satisfied for $k+1 \le i < j \le k+l$. That
each $\alpha_{e_i}$ is an automorphism guarantees that
Condition~(1) is satisfied for $1 \le i \le k$ and $1 \le j \le
l$. Condition~(1) is satisfied for $1 \le i < j \le k$ because
Condition~(1) is satisfied in $E_\Lambda$. Associativity of
composition in $\Lambda$ and that $\alpha$ is an action ensure
that the above factorisation rules satisfy Condition~(2).

It follows from \cite[Example~1.5(4) and Theorems 2.1~and~2.2]{FS1}
that there is a unique $(k+l)$-graph $\cpg{\Lambda}{\alpha}{l}$ with
skeleton $(E,F)$. Fix $(p,m) \in \NN^{k+l}$, and $\lambda \in
\Lambda^p$. Factorise $\lambda = f_1 \dots f_{|p|}$ as a sequence of
edges from $E_\Lambda$. Let $|m|$ denote the \emph{length} $m_1 +
\dots + m_k$ of $m$ as an element of $\NN^l$ with respect to the
usual basis, and fix $a_1, \dots, a_{|m|} \in \{1, \dots, l\}$ such
that $m = e_{a_1} + e_{a_2} + \dots + e_{a_{|m|}}$. This gives us a
path
\[
f_1 \dots f_{|p|} (s(\lambda), a_1) (\alpha_{a_1}^{-1}(s(\lambda)), a_2)
\dots (\alpha_{m - a_{|m|}}^{-1}(s(\lambda)), a_{|m|})
\]
in $E$, and hence a path $\xi(\lambda,m) \in
(\cpg{\Lambda}{\alpha}{l})^{(p,m)}$. One checks using the
definition of $F$ that $\xi(\lambda,m)$ does not depend on the
choice of factorisation of $\lambda$ into edges or the
decomposition of $m$ into generators. Using the factorisation
property in $\cpg{\Lambda}{\alpha}{l}$, one checks that every
path in $(\cpg{\Lambda}{\alpha}{l})^{(p,m)}$ is of the form
$\xi(\lambda,m)$ for some $\lambda \in \Lambda^p$. It follows
that $\xi$ is a bijection between $\Lambda^p \times \{m\}$ and
$(\cpg{\Lambda}{\alpha}{l})^{(p,m)}$. One then checks using the
definition of $(E,F)$ that this bijection satisfies
(2)~and~(3).

To prove the final statement of the Proposition, we first claim
that for $\mu,\nu \in \Lambda$ and $m,n \in \NN^l$,
\begin{equation}\label{eq:MCE relationship}
    \MCE_{\cpg{\Lambda}{\alpha}{l}}((\mu,m),(\nu,n)) =
\MCE_\Lambda(\mu,\nu) \times \{m \vee n\}.
\end{equation}
To see this, suppose first that $(\lambda,p) \in
\MCE_{\cpg{\Lambda}{\alpha}{l}}((\mu,m), (\nu,n))$. Then
$d(\lambda,p) = (d(\mu)\vee d(\nu), m \vee n)$, which gives $p
= m \vee n$ and $d(\lambda) = d(\mu) \vee d(\nu)$. Moreover,
$(\lambda,p) = (\mu,m)(\eta,p-m)$ where $\mu\alpha_{m}(\eta) =
\lambda$. Likewise, $(\lambda,p) = (\nu,n)(\xi,p-n)$ where
$\nu\alpha_{n}(\xi) = \lambda$. Hence $\lambda \in
\MCE_\Lambda(\mu,\nu)$ and $(\lambda,p) \in
\MCE_\Lambda(\mu,\nu) \times \{m \vee n\}$.

Now suppose that $\lambda \in \MCE(\mu,\nu)$, and let $p := m
\vee n$. By definition, we have $d(\lambda,p) = d(\mu,m) \vee
d(\nu,n)$. Write $\lambda = \mu\eta = \nu\xi$. By definition of
composition in $\cpg{\Lambda}{\alpha}{l}$ we have $(\lambda,p)
= (\mu,m)(\alpha_m^{-1}(\eta), p-m) =
(\nu,n)(\alpha_n^{-1}(\xi), p-n)$, so $(\lambda,p) \in
\MCE_{\cpg{\Lambda}{\alpha}{l}}((\mu,m),(\nu,n))$. This
establishes~\eqref{eq:MCE relationship}, and in particular
implies that $\cpg{\Lambda}{\alpha}{l}$ is finitely aligned if
and only if $\Lambda$ is finitely aligned.

Next note that by construction of $\cpg{\Lambda}{\alpha}{l}$,
\[
|v\Lambda^p| = |(v,\Z{l})(\cpg{\Lambda}{\alpha}{l})^{(p,m)}|
\quad\text{for all $v \in \Lambda^0$ and $(p,m) \in \NN^{k+l}$}.
\]
In particular, $\cpg{\Lambda}{\alpha}{l}$ is row-finite if and
only if $\Lambda$ is row-finite, and $\cpg{\Lambda}{\alpha}{l}$
has no sources if and only if $\Lambda$ has no sources.
\end{proof}

\begin{examples}
\begin{enumerate}
\item Let $k = l = 1$, and let $\Lambda$ be the $1$-graph
    consisting of a two-sided infinite path with vertices
    $\{v_n : n \in \ZZ\}$ and edges $\{f_n : n \in \ZZ\}$
    where $r(f_n) = v_n$ and $s(f_n) = v_{n+1}$. Let
    $\alpha$ be the automorphism of $\Lambda$ satisfying
    $\alpha(v_n) = v_{n+2}$ and $\alpha(f_n) = f_{n+2}$.
    Then the skeleton of $\cpg{\Lambda}{\alpha}{l}$ is as
    follows:
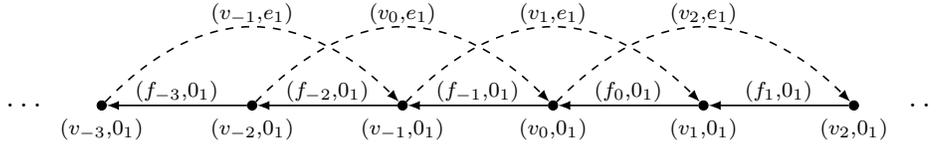
\[
\begin{tikzpicture}[scale=2]
    %define the vertex command - places a bullet with "anchors"
    %labelled #1 at position (#2,#3)
    \def\vertex(#1) at (#2,#3){
        \node[inner sep=0pt, circle, fill=black] (#1) at (#2,#3)
        [draw] {.};
    }
    %place vertices with labels underneath
    \vertex(vertA) at (-3,0)
    \node[inner sep=3pt, anchor = north] at (vertA.south)
    {$\scriptstyle (v_{-3},\Z{1})$};
    \vertex(vertB) at (-2,0)
    \node[inner sep=3pt, anchor = north] at (vertB.south)
    {$\scriptstyle (v_{-2},\Z{1})$};
    \vertex(vertC) at (-1,0)
    \node[inner sep=3pt, anchor = north] at (vertC.south)
    {$\scriptstyle (v_{-1},\Z{1})$};
    \vertex(vertD) at (0,0)
    \node[inner sep=3pt, anchor = north] at (vertD.south)
    {$\scriptstyle (v_{0},\Z{1})$};
    \vertex(vertE) at (1,0)
    \node[inner sep=3pt, anchor = north] at (vertE.south)
    {$\scriptstyle (v_{1},\Z{1})$};
    \vertex(vertF) at (2,0)
    \node[inner sep=3pt, anchor = north] at (vertF.south)
    {$\scriptstyle (v_{2},\Z{1})$};
%    \vertex(vertG) at (3,0)
%    \node[inner sep=3pt, anchor = north] at (vertG.south)
%    {$\scriptstyle (v_{3},\Z{1})$};
    %draw solid edges with labels above
    \draw[style=semithick,-latex] (vertB.west)--(vertA.east);
    \node[inner sep=1pt, anchor = south] at (-2.5,0) {$\scriptstyle
    (f_{-3}, \Z{1})$};
    \draw[style=semithick,-latex] (vertC.west)--(vertB.east);
    \node[inner sep=1pt, anchor = south] at (-1.5,0) {$\scriptstyle
    (f_{-2}, \Z{1})$};
    \draw[style=semithick,-latex] (vertD.west)--(vertC.east);
    \node[inner sep=1pt, anchor = south] at (-0.5,0) {$\scriptstyle
    (f_{-1}, \Z{1})$};
    \draw[style=semithick,-latex] (vertE.west)--(vertD.east);
    \node[inner sep=1pt, anchor = south] at (0.5,0) {$\scriptstyle
    (f_{0}, \Z{1})$};
    \draw[style=semithick,-latex] (vertF.west)--(vertE.east);
    \node[inner sep=1pt, anchor = south] at (1.5,0) {$\scriptstyle
    (f_{1}, \Z{1})$};
%    \draw[style=semithick,-latex] (vertG.west)--(vertF.east);
%    \node[inner sep=1pt, anchor = south] at (2.5,0) {$\scriptstyle
%    (f_{2}, \Z{1})$};
    %draw dashed edges with labels above
    \draw[style=semithick, style=dashed, -latex] (vertA.north east)
    parabola[parabola height=0.5cm] (vertC.north west);
    \node[inner sep=2pt, anchor = south] at (-2,0.5) {$\scriptstyle
    (v_{-1}, e_1)$};
    \draw[style=semithick, style=dashed, -latex] (vertB.north east)
    parabola[parabola height=0.5cm] (vertD.north west);
    \node[inner sep=2pt, anchor = south] at (-1,0.5) {$\scriptstyle
    (v_{0}, e_1)$};
    \draw[style=semithick, style=dashed, -latex] (vertC.north east)
    parabola[parabola height=0.5cm] (vertE.north west);
    \node[inner sep=2pt, anchor = south] at (0,0.5) {$\scriptstyle
    (v_{1}, e_1)$};
    \draw[style=semithick, style=dashed, -latex] (vertD.north east)
    parabola[parabola height=0.5cm] (vertF.north west);
    \node[inner sep=2pt, anchor = south] at (1,0.5) {$\scriptstyle
    (v_{2}, e_1)$};
%    \draw[style=semithick, style=dashed, -latex] (vertE.north east)
%    parabola[parabola height=0.5cm] (vertG.north west);
%    \node[inner sep=2pt, anchor = south] at (2,0.5) {$\scriptstyle
%    (v_{3}, e_1)$};
    %place dots at ends
    \node at (-3.5,0) {\dots};
    \node at (2.5,0) {\dots};
\end{tikzpicture}
\]
\item Fix $\theta \in [0,1] \setminus \QQ$, and let $[c_1,
    c_2, \dots]$ be its reduced continued fraction
    expansion. For each $n \in \NN$, let $\Phi_n$ be the
    matrix $\big(\begin{smallmatrix} c_n & 1 \\ 1 & 0
    \end{smallmatrix}\big)$. Let $T_n$ be the sequence of triangular
    numbers $T_n := \sum_{i=1}^n i = \frac{n(n+1)}{2}$, and
    for $n \in \NN$ let $A_n$ denote the matrix
    $\prod^{T_{n}}_{i = T_{n-1} + 1} \Phi_i$. Let
    $a^n_{i,j}$ denote the $i\,j$-entry of $A_n$ for $i,j =
    1,2$. Let $E$ be the directed graph with vertices
    $\{v^m_i : m \in \NN, i \in \{1,2\}\}$, and with edges
    $\{e^m_{i,j}(n) : m \in \NN, i,j \in \{1,2\}, n \in
    \ZZ/a^n_{i,j}\ZZ\}$, where $r(e^m_{i,j}(n) = v^m_i$ and
    $s(e^m_{i,j}(n)) = v^{m+1}_j$. So $E$ consists of the
    vertices and solid edges in the diagram below (where a
    label $n$ on an arrow indicates a bundle of $n$
    parallel edges).
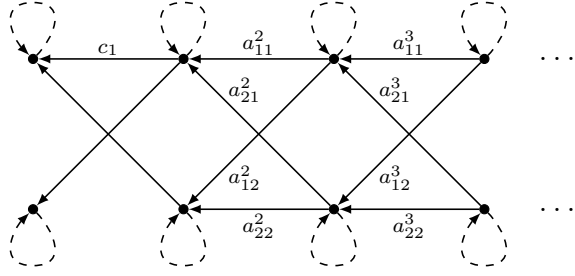
\[
\begin{tikzpicture}
    \def\vertex(#1) at (#2,#3){
        \node[inner sep=0pt, circle, fill=black] (#1) at (#2,#3)
        [draw] {.};
    }
    \vertex(11) at (0, 1)
    \vertex(12) at (0, -1)
    \vertex(21) at (2, 1)
    \vertex(22) at (2, -1)
    \vertex(31) at (4, 1)
    \vertex(32) at (4, -1)
    \vertex(41) at (6, 1)
    \vertex(42) at (6, -1)
    \node at (7, 1) {$\dots$};
    \node at (7, -1) {$\dots$};
    \draw[style=semithick, -latex] (21.west)--(11.east) node[pos=0.5, anchor=south, inner sep=1pt]{$\scriptstyle c_1$};
    \draw[style=semithick, -latex] (21.south west)--(12.north east);
    \draw[style=semithick, -latex] (22.north west)--(11.south east);
%    \draw[style=semithick, -latex] (22.west)--(12.east);
%
    \draw[style=semithick, -latex] (31.west)--(21.east) node[pos=0.5, anchor=south, inner sep=1pt]{$\scriptstyle a^2_{11}$};
    \draw[style=semithick, -latex] (31.south west)--(22.north east) node[pos=0.7, circle, anchor=north west, inner sep=0pt]{$\scriptstyle a^2_{12}$};
    \draw[style=semithick, -latex] (32.north west)--(21.south east) node[pos=0.7, circle, anchor=south west, inner sep=0pt]{$\scriptstyle a^2_{21}$};
    \draw[style=semithick, -latex] (32.west)--(22.east) node[pos=0.5, anchor=north, inner sep=1.5pt]{$\scriptstyle a^2_{22}$};
    \draw[style=semithick, -latex] (41.west)--(31.east) node[pos=0.5, anchor=south, inner sep=1pt]{$\scriptstyle a^3_{11}$};
    \draw[style=semithick, -latex] (41.south west)--(32.north east) node[pos=0.7, circle, anchor=north west, inner sep=0pt]{$\scriptstyle a^3_{12}$};
    \draw[style=semithick, -latex] (42.north west)--(31.south east) node[pos=0.7, circle, anchor=south west, inner sep=0pt]{$\scriptstyle a^3_{21}$};
    \draw[style=semithick, -latex] (42.west)--(32.east) node[pos=0.5, anchor=north, inner sep=1.5pt]{$\scriptstyle a^3_{22}$};
    \draw[style=semithick, style=dashed, -latex] (11.north east)
        .. controls (0.25,1.25) and (0.5,1.75) ..
        (0,1.75)
        .. controls (-0.5,1.75) and (-0.25,1.25) ..
        (11.north west);
    \draw[style=semithick, style=dashed, -latex] (21.north east)
        .. controls (2.25,1.25) and (2.5,1.75) ..
        (2,1.75)
        .. controls (1.5,1.75) and (1.75,1.25) ..
        (21.north west);
    \draw[style=semithick, style=dashed, -latex] (31.north east)
        .. controls (4.25,1.25) and (4.5,1.75) ..
        (4,1.75)
        .. controls (3.5,1.75) and (3.75,1.25) ..
        (31.north west);
    \draw[style=semithick, style=dashed, -latex] (41.north east)
        .. controls (6.25,1.25) and (6.5,1.75) ..
        (6,1.75)
        .. controls (5.5,1.75) and (5.75,1.25) ..
        (41.north west);
    \draw[style=semithick, style=dashed, -latex] (12.south east)
        .. controls (0.25,-1.25) and (0.5,-1.75) ..
        (0,-1.75)
        .. controls (-0.5,-1.75) and (-0.25,-1.25) ..
        (12.south west);
    \draw[style=semithick, style=dashed, -latex] (22.south east)
        .. controls (2.25,-1.25) and (2.5,-1.75) ..
        (2,-1.75)
        .. controls (1.5,-1.75) and (1.75,-1.25) ..
        (22.south west);
    \draw[style=semithick, style=dashed, -latex] (32.south east)
        .. controls (4.25,-1.25) and (4.5,-1.75) ..
        (4,-1.75)
        .. controls (3.5,-1.75) and (3.75,-1.25) ..
        (32.south west);
    \draw[style=semithick, style=dashed, -latex] (42.south east)
        .. controls (6.25,-1.25) and (6.5,-1.75) ..
        (6,-1.75)
        .. controls (5.5,-1.75) and (5.75,-1.25) ..
        (42.south west);
\end{tikzpicture}
\]
Define an automorphism $\alpha_1$ of $E$ as follows:
$\alpha_1(v^m_i) = v^m_i$ for all $v \in E^0$, and
$\alpha_1(e^m_{i\,j}(n)) := e^m_{i\,j}(n + 1)$. That is,
$\alpha_1$ fixes all the vertices, and cyclicly permutes
parallel edges. Let $\Lambda = E^*$ be the path-category of
$E$ regarded as a $1$-graph. Then $\alpha_1$ extends
uniquely to an automorphism $\bar{\alpha}_1$ of $\Lambda$.
Let $\bar\alpha$ be the action of $\ZZ$ on $\Lambda$
generated by $\bar\alpha_1$. The skeleton of
$\cpg{\Lambda}{\bar\alpha}{}$ is the $2$-coloured graph
pictured above, which is identical to the one in
\cite[Figure~3]{PRRS}. In particular, it follows from
\cite[Example~6.5]{PRRS} that the $C^*$-algebra of this
$2$-graph is Morita equivalent to the irrational rotation
algebra $A_\theta$.
\end{enumerate}
\end{examples}

\begin{theorem} \label{thm:crossed product}
Let $\Lambda$ be a finitely aligned $k$-graph, and $\alpha$ an
action of $\ZZ^l$ by automorphisms of $\Lambda$. Let
$\tilde\alpha$ be the corresponding action of $\ZZ^l$ on
$C^*(\Lambda)$ as in Proposition~\ref{prp:obtain action}, and
denote by $\pi : C^*(\Lambda) \to C^*(\Lambda)
\times_{\tilde\alpha} \ZZ^l$ and by $U : \ZZ^l \to
\Mm(C^*(\Lambda) \times_{\tilde\alpha} \ZZ^l)$ the universal
generating covariant representation of the dynamical system
$(C^*(\Lambda), \ZZ^l, \tilde\alpha)$. There is a unique
isomorphism $\phi : C^*(\cpg{\Lambda}{\alpha}{l}) \to
C^*(\Lambda) \times_{\tilde\alpha} \ZZ^l$ which satisfies
$\phi(s_{(\lambda,m)}) := \pi(s_\lambda) U_m$ for each $\lambda
\in \Lambda$ and $m \in \NN^l$.
\end{theorem}
\begin{proof}
For $(\lambda,m) \in \cpg{\Lambda}{\alpha}{l}$, let $t_{(\lambda,m)}
:= \pi(s_\lambda) U_m$. We claim that $\{t_{(\lambda,m)} :
(\lambda,m) \in \cpg{\Lambda}{\alpha}{l}\}$ is a Cuntz-Krieger
$\cpg{\Lambda}{\alpha}{l}$-family. First note that $U_{\Z{l}}$ is the
identity element of $\Mm(C^*(\Lambda) \times_{\tilde\alpha} \ZZ^l)$,
so $t_{(\lambda,\Z{l})} = \pi(s_\lambda)$ for each $\lambda \in
\Lambda$. Since $\pi$ is an isomorphism, it follows that
$\{t_{(\lambda,\Z{l})} : \lambda \in \Lambda\}$ forms a Cuntz-Krieger
$\Lambda$-family in $C^*(\Lambda) \times_\alpha \ZZ^l$. In
particular, the elements $\{t_{(v,\Z{l})} : (v,\Z{l}) \in
(\cpg{\Lambda}{\alpha}{l})^{\Z{k+l}}\}$ are mutually orthogonal
projections, which establishes~(TCK1). For $(\mu,m)$ and $(\nu,n)$
with $r(\nu) = \alpha_m^{-1}(s(\mu))$, we have
\begin{align*}
    t_{(\mu,m)}t_{(\nu,n)}
        &= \pi(s_\mu)U_m\pi(s_\nu)U_n \\
        &= \pi(s_\mu)U_m\pi(s_\nu)U_m^* U_mU_n \\
        &= \pi(s_\mu)\pi(\tilde\alpha_m(s_\nu)) U_m U_n  \\
        &= \pi(s_{\mu\alpha_m(\nu)}) U_{m+n} \\
        &= t_{(\mu,m)(\nu,n)}.
\end{align*}
This establishes~(TCK2).

To show that~(TCK3) holds, fix $(\mu,m), (\nu,n) \in
\cpg{\Lambda}{\alpha}{l}$ with $r(\mu,m) = r(\nu,n)$. We
calculate:
\begin{align*}
    t_{(\mu,m)} t^*_{(\mu,m)} t_{(\nu,n)} t^*_{(\nu,n)}
        &= \pi(s_\mu)U_m U^*_m \pi(s_\mu)^*\pi(s_\nu)U_n U^*_n
        \pi(s_\nu)^*  \\
        &= \pi\Big(\sum_{\lambda \in \MCE(\mu,\nu)} s_\lambda
        s^*_\lambda\Big)
\end{align*}
because $U_m, U_n$ are unitaries and $\{\pi(s_\lambda) :
\lambda \in \Lambda\}$ is a Cuntz-Krieger $\Lambda$-family.
Hence~\eqref{eq:MCE relationship} implies that
\[
    t_{(\mu,m)} t^*_{(\mu,m)} t_{(\nu,n)} t^*_{(\nu,n)} =
    \sum_{(\lambda,p) \in
        \MCE_{\cpg{\Lambda}{\alpha}{l}}((\mu,m),(\nu,n))}
        t_{(\lambda,p)} t^*_{(\lambda,p)},
\]
and multiplying both sides of this equation on the left by
$t^*_{(\mu,m)}$ and on the right by $t_{(\nu,n)}$ gives~(TCK3).

To show that~(CK) holds, we first establish that finite
exhaustive sets in $\cpg{\Lambda}{\alpha}{l}$ project onto
finite exhaustive sets in $\Lambda$. Fix a vertex $(v,\Z{l})$
of $\cpg{\Lambda}{\alpha}{l}$, and let $F$ be an exhaustive
subset of $(v,\Z{l})\Lambda^{\alpha}$. Let $F_1 := \{\lambda
\in \Lambda : \text{ there exists } m \in \NN^l \text{ such
that } (\lambda,m) \in F\}$, so that $F_1$ is the projection of
$F$ onto $\Lambda$. We claim that $F_1$ is exhaustive in
$\Lambda$. To see this, fix $\mu \in \Lambda$ with $r(\mu) =
v$. Then $(\mu,\Z{l}) \in (v,\Z{l})(\cpg{\Lambda}{\alpha}{l})$.
Hence there exists $(\lambda,p) \in F$ with
$\MCE_{\cpg{\Lambda}{\alpha}{l}}((\lambda,p),(\mu,\Z{l})) \not=
\emptyset$. Equation~\eqref{eq:MCE relationship} implies that
$\MCE_{\cpg{\Lambda}{\alpha}{l}}((\lambda,p),(\mu,\Z{l})) =
\MCE_\Lambda(\lambda,\mu) \times \{p\}$, so
$\MCE_\Lambda(\lambda,\mu) \not= \emptyset$. Since $\mu \in
v\Lambda$ was arbitrary, and since $\lambda \in F_1$, it
follows that $F_1$ is exhaustive. Now to establish~(CK),
suppose that $F$ is finite exhaustive in
$\cpg{\Lambda}{\alpha}{l}$, so that, by the above, $F_1$ is
finite exhaustive in $\Lambda$. Then
\begin{align*}
    \prod_{(\lambda,p) \in F}(t_{(v,\Z{l})} - t_{(\lambda,p)}
    t^*_{(\lambda,p)})
        &= \prod_{(\lambda,p) \in F}(\pi(s_v)U_0 - \pi(s_\lambda)
        U_p
        U^*_p \pi(s_\lambda)^*) \\
        &= \pi\Big(\prod_{\lambda \in F_1} (s_v - s_\lambda
        s^*_\lambda)\Big) \\
        &= 0
\end{align*}
because $s_\Lambda$ satisfies relation~(CK).

We have now proved that $t_{\cpg{\Lambda}{\alpha}{l}}$ is a
Cuntz-Krieger $(\cpg{\Lambda}{\alpha}{l})$-family.

The universal property of $C^*(\cpg{\Lambda}{\alpha}{l})$
implies that there is a homomorphism $\phi :
C^*(\cpg{\Lambda}{\alpha}{l}) \to C^*(\Lambda)
\times_{\tilde\alpha} \ZZ^l$ which satisfies
$\phi(s_{(\lambda,m)}) = \pi(s_\lambda)U_m$ for all $\lambda
\in \Lambda$, $m \in \NN^l$.

We claim that $\phi$ is surjective. The crossed product
$C^*(\Lambda) \times_{\tilde\alpha} \ZZ^l$ is, by definition,
generated by elements of the form $\pi(a)U_m$ where $a \in
C^*(\Lambda)$ and $m \in \ZZ^k$. Hence it suffices to show that
$\pi(a)U_m$ is in the image of $\phi$ for each $a \in
C^*(\Lambda)$ and $m \in \ZZ^l$. Since $C^*(\Lambda) =
\clsp\{s_\mu s^*_\nu : \mu,\nu \in \Lambda$, it therefore
suffices to show that $s_\mu s^*_\nu U_m$ is in the range of
$\phi$ for each $m \in \ZZ^l$ and $\mu,\nu \in \Lambda$. Fix
$\mu,\nu \in \Lambda$ and $m \in \ZZ^l$. Write $m = m_+ - m_-$
where $m_+, m_- \in \NN^l$. Since $U^*_m s^*_\nu U_m =
\tilde{\alpha}_{-m}(s^*_\nu) = s^*_{\alpha_{-m}(\nu)}$, we have
\[
s_\mu s^*_\nu U_m
= s_\mu U_m U^*_m s^*_\nu U_m
= s_{\mu} U_{m_+} U^*_{m_-} s^*_{\alpha_{-m}(\nu)}
%= \big(s_{\mu} U_{m_+}\big)\big(s_{\alpha_{-m}(\nu)}U_{m_-}\big)^* \\
= \phi\big(t_{(\mu,m_+)} t_{(\alpha_{-m}(\nu), m_-)}^*\big).
\]
Hence $\phi$ is surjective. It remains only to show that $\phi$
is injective.

Let $\gamma$ denote the gauge action of $\TT^k$ on
$C^*(\Lambda)$ and let $\gamma^\alpha$ denote the gauge action
of $\TT^{k+l}$ on $C^*(\cpg{\Lambda}{\alpha}{l})$. The
universal property of the crossed product $C^*(\Lambda)
\times_{\tilde\alpha} \ZZ^l$ can be used to deduce that there
is an action $\bar\gamma$ of $\TT^k$ on $C^*(\Lambda)
\times_{\tilde\alpha} \ZZ^l$ which satisfies
$\overline{\gamma_z}(\pi(a)U_m) = \pi(\gamma_z(a))U_m$ for all
$z \in \TT^k$, $a \in C^*(\Lambda)$, and $m \in \ZZ^l$. Let
$\widehat{\tilde{\alpha}}$ denote the dual action of $\TT^l =
\widehat{\ZZ^l}$ on $C^*(\Lambda) \times_{\tilde\alpha} \ZZ^l$
which satisfies $\widehat{\tilde{\alpha}}_w(\pi(a)U_m) :=
w^m\pi(a)U_m$. Identifying $\TT^{k+l}$ with $\{(z,w) : z \in
\TT^k, w \in \TT^l\}$, define automorphisms $\{\beta_{(z,w)} :
(z,w) \in \TT^{k+l}\}$ of $C^*(\Lambda) \times_{\tilde\alpha}
\ZZ^l$ by $\beta_{(z,w)} := \overline{\gamma_z} \circ
\widehat{\tilde{\alpha}}_w$. It is easy to see that $\beta$
determines an action of $\TT^{k+l}$ on the crossed product
algebra, and one can check on generators that $\phi \circ
\gamma^\alpha_{(z,w)} = \beta_{(z,w)} \circ \phi$ for all
$(z,w) \in \TT^{k+l}$. Since $\pi$ is injective, we have $t_v =
\pi(s_v) \not= 0$ for all $v$. The gauge-invariant uniqueness
theorem~\cite[Theorem~4.2]{RSY2} therefore implies that $\phi$
is injective.
\end{proof}

\subsection{Recognising crossed-product $k$-graphs}
There is a converse of sorts to
Proposition~\ref{prp:Lambda^alpha}; that is, one can tell by
looking at a $(k+l)$-graph whether or not it is of the form
$\cpg{\Lambda}{\alpha}{l}$ for some $k$-graph $\Lambda$ and
some action $\alpha$ of $\ZZ^l$ on $\Lambda$ by automorphisms.

We require some notation. Given a $(k+l)$-graph $\Xi$, we write
$\Xi^{(\NN^k,\Z{l})}$ for the $k$-graph with morphisms
$\bigcup_{p \in \NN^k} \Xi^{(p,\Z{l})}$ and degree functor
$d^{(\NN^k,\Z{l})}(\xi) = (d(\xi)_1, \dots, d(\xi)_k) \in
\NN^k$.

\begin{prop}\label{prp:when cpgraph}
Let $\Xi$ be a $(k+l)$-graph. Suppose that for every $v \in
\Xi^{\Z{k+l}}$ and every $j \in \{1, \dots, l\}$, we have
$|v\Xi^{(\Z{k}, e_j)}| = |\Xi^{(\Z{k}, e_j)}v| = 1$. Then for each
vertex $v \in \Xi^{\Z{k+l}}$ and $m \in \NN^l$ there is a unique path
$\eta_{v,m}$ in $\Xi^{(\Z{k},m)} v$, and there is a unique action
$\alpha$ of $\ZZ^l$ on $\Xi^{(\NN^k,\Z{l})}$ satisfying
$\alpha_m(\xi) := (\eta_{r(\xi),m} \xi)(\Z{k+l}, d(\xi))$ for all
$\xi \in \Xi^{(\NN^k, \Z{l})}$ and $m \in \NN^l$. Moreover, $\Xi$ is
isomorphic to $\cpg{\Xi^{(\NN^k,\Z{l})}}{\alpha}{l}$.
\end{prop}
\begin{proof}
Since $|\Xi^{(\Z{k}, e_j)}v| = 1$ for all $v$ and $j$, it is
immediate from the factorisation property that
$|\Xi^{(\Z{k},m)} v| = 1$ for all $m \in \NN^l$ and $v \in
\Xi^{\Z{k+l}}$. Arguments like those of \cite[Lemma~3.3]{PRRS}
show that for $m \in \NN^l$, the formula
\[
\xi \mapsto (\eta_{r(\xi),m} \xi)(\Z{k+l}, d(\xi))
\]
determines an automorphism $\alpha_m$ of $\Xi^{(\NN^k,\Z{l})}$
for each $m \in \NN^l$ and that $\alpha_m \circ \alpha_n =
\alpha_{m+n}$ for all $m,n \in \NN^l$. Hence the $\alpha_m$
determine an action $\alpha$ of $\ZZ^l$ on
$\Xi^{(\NN^k,\Z{l})}$ as claimed: writing $m \in \ZZ^l$ as $m =
m_+ - m_-$ where $m_+, m_- \in \NN^l$, we define $\alpha_{m} :=
\alpha_{m_+} \circ \alpha^{-1}_{m_-}$. It is easy to check that
the skeleton of $\cpg{\Xi^{(\NN^k,\Z{l})}}{\alpha}{l}$ is
identical to that of $\Xi$, so the isomorphism $\Xi \cong
\cpg{\Xi^{(\NN^k,\Z{l})}}{\alpha}{l}$ follows from the
uniqueness assertion of \cite[Theorem~2.1]{FS1}.
\end{proof}

\begin{examples}\label{egs:recognise cpgs}
(1) As in Section~3 of \cite{PRRS}, let $\Lambda$ be a row-finite
$2$-graph with no sources such that $\Lambda^{(\NN,0)}$ contains no
cycles and each vertex $v \in \Lambda^{\Z{2}}$ is the range of an
isolated cycle in $\Lambda^{(0,\NN)}$. It was shown in
\cite[Theorem~3.1]{PRRS}, $C^*(\Lambda)$ is an A$\TT$-algebra.

By Proposition~\ref{prp:when cpgraph}, we see that $C^*(\Lambda)$ is
isomorphic to a crossed product by $\ZZ$ of the AF algebra
$C^*(\Lambda^{(\NN,0)})$. In Section~\ref{sec:K-theory} we show how
to use the Pimsner-Voiculescu exact sequence to calculate the
$K$-theory of $C^*(\Lambda)$.

(2) Consider the $3$-graphs $\tgrphlim(\Delta_2/H_n, p_n)$
discussed in \cite[Section~6.4]{KPS}. Since each $\Delta_2/H_n$
is a quotient of $\Delta_2$, each vertex in
$\tgrphlim(\Delta_2/H_n, p_n)$ both emits and receives exactly
one edge of degree $e_1$ and exactly one edge of degree $e_2$.
After a change of basis for $\NN^3$, we can therefore use
Proposition~\ref{prp:when cpgraph} and Theorem~\ref{thm:crossed
product} to realise $C^*(\tgrphlim(\Delta_2/H_n, p_n))$ as a
crossed product of the AF algebra $C^*(\tgrphlim(\Delta_2/H_n,
p_n)^{(\Z{2},\NN)})$ by $\ZZ^2$. Indeed, taking the corner
generated by the lone vertex in $\Delta_2/H_1$, we recover the
crossed product of $C_0(\varprojlim \ZZ^2/H_n)$ by a
generalised odometer action discussed in
\cite[Remark~6.12]{KPS}.

(3) Fix an integer $l \ge 1$. It is easy to verify that there
is an isomorphism of $C^*(\Delta_l)$ onto $\Kk(\ell^2(\ZZ^l))$
which takes $s_{(m, m+n)}$ to the matrix unit $\theta_{m,
m+n}$; we will henceforth identify $C^*(\Delta_l)$ with
$\Kk(\ell^2(\ZZ^l))$ via this isomorphism.

Fix a $k$-graph $\Lambda$. Consider the Cartesian product
$(k+l)$-graph $\Lambda \times \Delta_l := \{(\lambda, (m,n)) :
\lambda \in \Lambda, (m,n) \in \Delta_l\}$ with coordinatewise
range, source and composition maps, and degree map
$d(\lambda,(m,n)) = (d(\lambda),d(m,n)) = (d(\lambda),n-m)$
(see \cite[Proposition~1.8]{KP}). Clearly each vertex of
$\Lambda \times \Delta_l$ emits and receives exactly one edge
of degree $e_j$ for $k+1 \le j \le k+l$. Moreover,
$\Lambda^{(\NN^k, \Z{l})} \cong \bigsqcup_{m \in \ZZ^l} \Lambda
\times \{m\}$ is a disjoint union of copies of $\Lambda$
indexed by $\ZZ^l$, and the action $\alpha$ on this $k$-graph
arising from Proposition~\ref{prp:when cpgraph} is implemented
by translation in the $\ZZ^l$ coordinate; that is,
$\alpha_m'(\lambda, (m,n)) = (\lambda, (m+m', n+m'))$. Hence,
\[
C^*(\Lambda^{(\NN^k, \Z{l})}) \cong \bigoplus_{z \in \ZZ^l}
C^*(\Lambda) \cong C^*(\Lambda) \otimes c_0(\ZZ^l),
\]
and under this identification, $\tilde\alpha$ becomes $\id
\otimes \operatorname{lt}$. Since $c_0(\ZZ^l)
\times_{\operatorname{lt}} \ZZ^l$ is canonically isomorphic to
$\Kk(\ell^2(\ZZ^l))$, and since $C^*(\Delta_l)$ is also
canonically isomorphic to $\Kk(\ell^2(\ZZ^l))$
Theorem~\ref{thm:crossed product} re-proves the isomorphism
$C^*(\Lambda \times \Delta_l) \cong C^*(\Lambda) \otimes
C^*(\Delta_l)$ obtained from \cite[Corollary~3.5(iv)]{KP}.
\end{examples}

\subsection{Takai duality}\label{sec:takai} In this section we show how our
construction together with the skew-product construction of
\cite{KP, PQR2} provides a graph-theoretic realisation of Takai
duality for the $\ZZ^l$ actions discussed in
Theorem~\ref{thm:crossed product}.

We require some background regarding the skew-product of a
$k$-graph by an abelian group $G$ and its relationship to a
crossed product by an action of the dual group $\widehat{G}$
(see \cite[Section~5]{KP}). This construction has been
generalised to nonabelian groups $G$ using coactions (see
\cite{PQR2}), but for our purposes, the generality of \cite{KP}
suffices.

Let $\Xi$ be a $k$-graph, and let $c : \Xi \to G$ be a cocycle
into an abelian group $(G, +)$; that is, $c(\mu\nu) = c(\mu) +
c(\nu)$ whenever $s(\mu) = r(\nu)$. Following the conventions
of \cite{PQR2}, we define the \emph{skew-product $k$-graph}
$\Xi \times_c G$ to be equal as a set to $\Xi \times G$, with
structure maps
\begin{gather*}
s_c(\mu,g) := (s(\mu),g) \qquad r_c(\mu,g) = (r(\mu), c(\mu) + g) \\
(\mu,c(\nu) + g)(\nu,g) := (\mu\nu,g) \qquad\text{and}\qquad
d_c(\mu,g) := d(\mu)
\end{gather*}
for all $\mu,\nu \in \Xi$ such that $s(\mu) = r(\nu)$, and all $g \in
G$.

There is an action $\delta^c$ of $\widehat{G}$ on $C^*(\Xi)$
satisfying $\delta^c_\phi(s_\mu) := \phi(c(\mu)) s_\mu$ for all
$\phi \in \widehat{G}$ and $\mu \in \Xi$. Corollary~5.3 of
\cite{KP} states that $C^*(\Xi \times_c G) \cong C^*(\Xi)
\times_{\delta^c} \widehat{G}$.

Given an action $\alpha$ of $\ZZ^l$ on a $k$-graph $\Lambda$, the map
$c : \cpg{\Lambda}{\alpha}{l} \to \ZZ^l$ defined by $c(\lambda,m) :=
-m$ for all $\lambda \in \Lambda$ and $m \in \NN^l$ is a cocycle (see
Proposition~\ref{prp:Lambda^alpha}). We may therefore form the
skew-product $(k+l)$-graph $(\cpg{\Lambda}{\alpha}{l}) \times_c
\ZZ^l$.

Finally, recall from \cite[Proposition~1.8]{KP} that given a
$k$-graph $\Lambda$ and an $l$-graph $\Gamma$, the cartesian product
$\Lambda \times \Gamma$ becomes a $(k+l)$-graph with structure maps
and degree functor defined coordinatewise.

\begin{theorem}\label{thm:Takai duality}
Fix a row-finite $k$-graph $\Lambda$ with no sources and an integer
$l \ge 1$. Let $\alpha$ be an action of $\ZZ^l$ on $\Lambda$. Then
the formula
\[
\rho((\lambda,m),n) := (\alpha_{n-m}(\lambda), (n-m, n))
\]
determines an isomorphism of the skew-product graph $(\Lambda
\times_\alpha \ZZ^l) \times_c \ZZ^l$ onto the cartesian product
$\Lambda \times \Delta_l$.
\end{theorem}
\begin{proof}
To establish that $\rho$ is an isomorphism of $(k+l)$-graphs, first
observe that it is bijective and degree-preserving by definition. We
therefore need only show that it preserves range, source and
composition. We have
\[\begin{split}
r(\rho((\lambda,m),n))
 &= r(\alpha_{n-m}(\lambda), (n-m, n)) = (\alpha_{n-m}(r(\lambda)),
 (n-m,n-m)) \\
 &= \rho((r(\lambda),\Z{l}), n-m) = \rho(r_c((\lambda, m),n))
\end{split}\]
and
\[\begin{split}
s(\rho((\lambda,m),n))
 &= s(\alpha_{n-m}(\lambda), (n-m, n)) = (\alpha_{n-m}(s(\lambda)),
 (n, n)) \\
 &= \rho((\alpha_{-m}(s(\lambda)),\Z{l}),n) =
 \rho(s_c((\lambda,m),n)),
\end{split}
\]
establishing that $\rho$ preserves the range and source maps. To see
that $\rho$ preserves composition, fix $\mu,\nu \in \Lambda$ and
$m,n,m' \in \NN^l$ such that $s(\mu) = \alpha_{m'}(r(\nu))$. Then
$((\mu,m'),n-m)$ and $((\nu, m),n)$ are composable in $(\Lambda
\times_\alpha \ZZ^l) \times_c \ZZ^l$ with
\[
((\mu,m'),n-m)((\nu, m),n) = ((\mu\alpha_{m'}(\nu), m+m'), n),
\]
and we must show that
\begin{equation} \label{eq:composition check}
\rho((\mu\alpha_{m'}(\nu), m+m'), n)
= \rho((\mu,m'),n-m) \rho((\nu, m),n).
\end{equation}
We calculate:
\begin{align*}
\rho((\mu\alpha_{m'}(\nu)&, m+m'), n)\\
 &= (\alpha_{n-m-m'}(\mu\alpha_{m'}(\nu)), (n-m-m', n)) \\
 &= (\alpha_{n-m-m'}(\mu)\alpha_{n-m}(\nu), (n-m-m', n-m)(n-m,n))\\
 &= (\alpha_{n-m-m'}(\mu), (n-m-m', n-m))(\alpha_{n-m}(\nu),(n-m,n))\\
 &= \rho((\mu, m'),n-m) \rho((\nu,m),n)
\end{align*}
as required.
\end{proof}

One interpretation of Theorem~\ref{thm:Takai duality} is that
for an action of $\ZZ^l$ on a $k$-graph $C^*$-algebra induced
by an action of $\ZZ^l$ on the $k$-graph itself, we may realise
Takai Duality at the level of higher-rank graphs.

\begin{cor}\label{cor:TD cor}
Let $\alpha$ be an action of $\ZZ^l$ on a row-finite $k$-graph
$\Lambda$ with no sources. Let $\widetilde{\alpha}$ denote the
induced action of $\ZZ^l$ on $C^*(\Lambda)$, and let
$\widehat{\widetilde{\alpha}}$ denote the dual action of $\TT^l$ on
$C^*(\Lambda) \times_{\widetilde{\alpha}} \ZZ^l$. The isomorphism
$\rho : (\cpg{\Lambda}{\alpha}{l}) \times_c \ZZ^l \to \Lambda \times
\Delta_l$ of Theorem~\ref{thm:Takai duality} induces an isomorphism
$\tilde\rho : (C^*(\Lambda) \times_{\widetilde{\alpha}} \ZZ^l)
\times_{\widehat{\widetilde{\alpha}}} \TT^l \to C^*(\Lambda) \otimes
\Kk(\ell^2(\ZZ^l))$.
\end{cor}
\begin{proof}
Corollary~3.5(iv) of \cite{KP} shows that $C^*(\Lambda \times
\Delta_l)$ is isomorphic to $C^*(\Lambda) \otimes C^*(\Delta_l)$, and
as mentioned in Examples~\ref{egs:recognise cpgs}(3), the map
$s_{(m,n)} \mapsto \theta_{m,n}$ determines an isomorphism of
$C^*(\Delta_l)$ onto $\Kk(\ell^2(\ZZ^l))$.
\end{proof}

\begin{rmk}
One can check that under appropriate conventions regarding dual
actions and crossed products, the isomorphism obtained from
Corollary~\ref{cor:TD cor} agrees with the Takai isomorphism as
described in, for example, \cite[Section~7.1]{TF^2B}.
\end{rmk}

\section{Simplicity of crossed products}\label{sec:simplicity}

In this section we investigate simplicity of
$C^*(\cpg{\Lambda}{\alpha}{l})$ when $\Lambda$ is row-finite
and has no sources.

Let $\Lambda$ be a row-finite $k$-graph with no sources. Recall
from \cite{KP} that $\Lambda$ is said to be \emph{cofinal} if
for every infinite path $x \in \Lambda^\infty$ and every vertex
$v \in \Lambda^0$ there is a vertex $x(n)$ on $x$ such that
$v\Lambda x(n) \not= \emptyset$. Recall from \cite{RobSi} that
$\Lambda$ has \emph{no local periodicity} if for every vertex
$v \in \Lambda^0$ and each pair of distinct elements $m,n \in
\NN^k$ there exists an infinite path $x \in v\Lambda^\infty$
such that $\sigma^m(x) \not= \sigma^n(x)$. Lemma~3.3 of
\cite{RobSi} implies that $\Lambda$ has no local periodicity if
and only if it satisfies the \emph{aperiodicity condition}
\cite[Condition~(A)]{KP}.

Recall that if $\phi \in \Aut(\Lambda)$ is an automorphism of a
$k$-graph $\Lambda$, then the formula
\[
\phi^\infty(x)(\Z{k},m) := \phi(x(\Z{k},m))\quad\text{for all $m
\in \NN^k$}
\]
defines a range-preserving bijection of $\Lambda^\infty$.

\begin{dfn}
Let $\Lambda$ be a row-finite $k$-graph with no sources and let
$\alpha$ be an action of $\ZZ^l$ on $\Lambda$ by automorphisms.
\begin{Enumerate}
\item We say that $\Lambda$ is \emph{$\alpha$-cofinal} if for every
    vertex $v \in \Lambda^{\Z{k}}$ and every infinite path $x \in
    \Lambda^\infty$ there exist $p \in \NN^k$ and $m,n \in \NN^l$
    such that $\alpha_{-m}(v)\Lambda \alpha_{-n}(x(p)) \not=
    \emptyset$. \item We say that $\Lambda$ is
    \emph{$\alpha$-aperiodic} if, for each vertex $v \in
    \Lambda^{\Z{k}}$ and each pair of distinct elements $(p,m)$ and
    $(q,n)$ of $\NN^k \times \NN^l$, there is a path $x \in
    v\Lambda^\infty$ such that $\sigma^p(\alpha^\infty_{-m}(x)) \not=
    \sigma^q(\alpha^\infty_{-n}(x))$.
\end{Enumerate}
\end{dfn}

\begin{theorem}\label{thm:first simplicity theorem}
Let $\Lambda$ be a row-finite $k$-graph with no sources, and let
$\alpha$ be an action of $\ZZ^l$ on $\Lambda$ by automorphisms. The
crossed-product $C^*$-algebra $C^*(\Lambda) \times_{\tilde\alpha}
\ZZ^k$ is simple if and only if $\Lambda$ is $\alpha$-cofinal and
$\alpha$-aperiodic.
\end{theorem}

To prove this theorem, we call upon the results of
\cite{RobSi}. We begin by showing how the infinite paths of
$\Lambda$ correspond to those of $\cpg{\Lambda}{\alpha}{l}$.
Specifically, we show that each infinite path in
$\cpg{\Lambda}{\alpha}{l}$ is determined by its restriction to
$\NN^k$.

\begin{lemma}\label{lem:(x, infty)->x}
Let $\Lambda$ be a row-finite $k$-graph with no sources, and suppose
that $\alpha$ is an action of $\ZZ^l$ on $\Lambda$ by automorphisms.
Suppose that $y,z \in (\cpg{\Lambda}{\alpha}{l})^\infty$ satisfy
$y((\Z{k},\Z{l}), (p,\Z{l})) = z((\Z{k},\Z{l}),(p,\Z{l}))$ for all $p
\in \NN^k$. Then $y = z$.
\end{lemma}
\begin{proof}
It suffices to show that $y((\Z{k},\Z{l}), (p,m)) = z((\Z{k},\Z{l}),
(p,m))$ for all $(p,m) \in \NN^k \times \NN^l$. To see this, fix
$(p,m) \in \NN^k \times \NN^l$. Since $y((\Z{k},\Z{l}),(p,\Z{l}))$
and $z((\Z{k},\Z{l}),(p,\Z{l}))$ coincide by assumption, they have
the same source. Since $(v,\Z{l})
(\cpg{\Lambda}{\alpha}{l})^{(\Z{k},m)}$ is a singleton set for any
fixed $(v,\Z{l}) \in (\cpg{\Lambda}{\alpha}{l})^{\Z{k+l}}$ and $m \in
\NN^l$, the paths $y((p,\Z{l}), (p,m))$ and $z((p,\Z{l}),(p,m))$ must
also coincide. Hence
\[\begin{split}
y((\Z{k},\Z{l}),(p,m)) &= y((\Z{k},\Z{l}), (p,\Z{l}))
y((p,\Z{l}),(p,m))\quad\text{and} \\
z((\Z{k},\Z{l}),(p,m)) &= (z(\Z{k},\Z{l}), (p,\Z{l}))
z((p,\Z{l}),(p,m))
\end{split}\]
are identical as required.
\end{proof}

Let $\Lambda$ be a row-finite $k$-graph with no sources. As in
\cite[Section~2]{KP}, the cylinder sets
$\lambda\Lambda^\infty$, $\lambda \in \Lambda$ form a basis of
compact open sets for a Hausdorff topology on $\Lambda^\infty$.

\begin{prop}\label{prp:inf path bijection}
Let $\Lambda$ be a row-finite $k$-graph with no sources, and
let $\alpha$ be an action of $\ZZ^l$ on $\Lambda$ by
automorphisms. Then there is a unique homeomorphism $x \mapsto
(x,\infty)$ from $\Lambda^\infty$ onto
$(\cpg{\Lambda}{\alpha}{l})^\infty$ such that
\begin{equation}\label{eq:(x,infty) def}
(x,\infty)((\Z{k},\Z{l}), (p,\Z{l})) = (x(\Z{k},p),
\Z{l})\quad\text{ for all $p \in\NN^k$.}
\end{equation}
\end{prop}
\begin{proof}
We first show that there exists a map $x \mapsto (x,\infty)$
satisfying~\eqref{eq:(x,infty) def}. To see this, fix $x \in
\Lambda^\infty$, and define paths $\{\lambda_{(p,m)} : (p,m) \in
\NN^{k+l}\} \subset \cpg{\Lambda}{\alpha}{l}$ by $\lambda_{(p,m)} :=
(x(\Z{k},p), m)$. If $(p,m) \le (q,n) \in \NN^{k+l}$, then we have
$\lambda_{(q,n)}((\Z{k},\Z{l}), (p,m)) = \lambda_{(p,m)}$, and it
follows from \cite[Remarks~2.2]{KP} that there is a unique infinite
path $(x,\infty) \in \cpg{\Lambda}{\alpha}{l}$ such that
$(x,\infty)((\Z{k},\Z{l}),(p,m)) = \lambda_{(p,m)}$ for all $(p,m)
\in \NN^{k+l}$. Since $\lambda_{(p,\Z{l})}$ is precisely the
right-hand side of~\eqref{eq:(x,infty) def}, this establishes the
existence of a the desired map $x \mapsto (x,\infty)$ from
$\Lambda^\infty$ to $(\cpg{\Lambda}{\alpha}{l})^\infty$.

Lemma~\ref{lem:(x, infty)->x} guarantees that $x \mapsto (x,\infty)$
is bijective and is the unique bijection
satisfying~\eqref{eq:(x,infty) def}.

It therefore remains only to show that $x \mapsto (x, \infty)$ is a
homeomorphism. To see this, observe that
$(\lambda,m)(\cpg{\Lambda}{\alpha}{l})^\infty =
(\lambda,\Z{l})(\cpg{\Lambda}{\alpha}{l})^\infty$ for all $\lambda
\in \Lambda$ and $m \in \NN^l$. In particular, the cylinder sets
$\{(\lambda,\Z{l})(\cpg{\Lambda}{\alpha}{l})^\infty : \lambda \in
\Lambda\}$ are a basis for the topology on
$(\cpg{\Lambda}{\alpha}{l})^\infty$, and since $x \mapsto (x,\infty)$
restricts to a bijection of $\lambda \Lambda^\infty$ onto
$(\lambda,\Z{l})(\cpg{\Lambda}{\alpha}{l})^\infty$, it follows that
$x \mapsto (x, \infty)$ is a homeomorphism.
\end{proof}

The next lemma shows how to express the shift maps on
$(\cpg{\Lambda}{\alpha}{l})^\infty$ in terms of the shift maps
on $\Lambda^\infty$ and the homeomorphisms $\alpha^\infty_p$ of
$\Lambda^\infty$ obtained from the automorphisms $\alpha_p$ of
$\Lambda$.

\begin{lemma}\label{lem:cp shift}
Let $\Lambda$ be a row-finite $k$-graph with no sources. and
let $\alpha$ be an action of $\ZZ^l$ on $\Lambda$ by
automorphisms. Then for $x \in \Lambda^\infty$ and $(p,m) \in
\NN^k \times \NN^l$, the shift map on
$(\cpg{\Lambda}{\alpha}{l})^\infty$ satisfies
$\sigma^{(p,m)}(x,\infty) = (\alpha^\infty_{-m}(\sigma^p(x)),
\infty)$. Moreover, $\sigma^p \circ \alpha^\infty_{-m} =
\alpha^\infty_{-m} \circ \sigma^p$; in particular,
$\sigma^{(p,m)}(x,\infty) = (\sigma^p(\alpha^\infty_{-m}(x)),
\infty)$.
\end{lemma}
\begin{proof}
Fix $(p,m) \in \NN^k \times \NN^l$ and $x \in \Lambda^\infty$.
Then $\sigma^{(p,m)}(x,\infty) =
\sigma^{(p,\Z{l})}(\sigma^{(\Z{k},m)}(x,\infty))$. For $q \in
\NN^k$, the initial segment of $(x,\infty)$ of degree $(q,m)$
is by definition equal to $(r(x), m)(\alpha_{-m}(x(\Z{k},q),
\Z{l}))$. Hence $\sigma^{(\Z{k},m)}(x,\infty) =
(\alpha^\infty_{-m}(x), \infty)$, and applying
$\sigma^{(p,\Z{l})}$ to both sides, we obtain the desired
identity $\sigma^{(p,m)}(x,\infty) =
(\alpha^\infty_{-m}(\sigma^p(x)), \infty)$. To see that
$\alpha^\infty_{-m} \circ \sigma^p = \sigma^p \circ
\alpha^\infty_{-m}$, fix $q \in \NN^k$ and calculate:
\[
\big(\alpha^\infty_{-m}(\sigma^p(x))\big)(0,q)
 = \alpha_{-m}(x(p,p+q))
 = \big(\alpha^\infty_{-m}(x)\big)(p,p+q)
 = \big(\sigma^p(\alpha^\infty_{-m}(x))\big)(0,q).
\]
This completes the proof.
\end{proof}

\begin{lemma}\label{lem:cofinality equivalence}
Let $\Lambda$ be a row-finite $k$-graph with no sources, and $\alpha$
an action of $\ZZ^l$ on $\Lambda$ by automorphisms. Then $\Lambda$ is
$\alpha$-cofinal if and only if $\cpg{\Lambda}{\alpha}{l}$ is cofinal
\end{lemma}
\begin{proof}
For $m \in \NN^l$ and $v,w \in \Lambda^{\Z{k}}$, we have $\lambda \in
\alpha_{-m}(v) \Lambda w$ if and only if $(\alpha_m(\lambda),m) \in
(v,\Z{l}) (\cpg{\Lambda}{\alpha}{l}) (w,\Z{l})$. Hence for $v, w \in
\Lambda^{\Z{k}}$, we have
\begin{equation} \label{eq:connectivity equivalence}
\text{$(v,\Z{l})(\cpg{\Lambda}{\alpha}{l}) (w,\Z{l}) \not=\emptyset$
if and only if $\alpha_{-m}(v)\Lambda w \not= \emptyset$ for some $m
\in \NN^l$.}
\end{equation}
Moreover, for $(p,n) \in \NN^{k + l}$, we have $\alpha_{-n}(x(p)) =
r(\alpha^\infty_{-n}(\sigma^p(x)))$, so
\begin{equation} \label{eq:vertex in inf path}
(\alpha_{-n}(x(p)),\Z{l}) = \sigma^{p,n}(x,\infty)(\Z{k},\Z{l}) =
(x,\infty)(p,n).
\end{equation}
Recall that every vertex $u$ of $\cpg{\Lambda}{\alpha}{l}$ is
of the form $(v,\Z{l})$ for some $v \in \Lambda^{\Z{k}}$.
Proposition~\ref{prp:inf path bijection} shows that every
infinite path $y$ of $\cpg{\Lambda}{\alpha}{l}$ is of the form
$(x,\infty)$ for some $x \in \Lambda^\infty$. Thus
\eqref{eq:connectivity equivalence}~and~\eqref{eq:vertex in inf
path} imply that $\Lambda$ is $\alpha$-cofinal if and only if,
for every vertex $u \in (\cpg{\Lambda}{\alpha}{l})^{\Z{k+l}}$
and every infinite path $y \in
(\cpg{\Lambda}{\alpha}{l})^\infty$, there exists $(p,n) \in
\NN^k \times \NN^l$ such that $u (\cpg{\Lambda}{\alpha}{l})
y(p,n) \not=\emptyset$, which is precisely the definition of
cofinality of $\cpg{\Lambda}{\alpha}{l}$.
\end{proof}

\begin{lemma}\label{lem:aperiodicity equivalence}
Let $\Lambda$ be a row-finite $k$-graph with no sources, and $\alpha$
an action of $\ZZ^l$ on $\Lambda$ by automorphisms. Then $\Lambda$ is
$\alpha$-aperiodic if and only if $\cpg{\Lambda}{\alpha}{l}$ has no
local periodicity in the sense of \cite{RobSi}.
\end{lemma}
\begin{proof}
The result follows from Lemma~\ref{lem:cp shift} and the definition
of $\alpha$-aperiodicity.
\end{proof}

\begin{rmk}
Suppose that the action $\alpha$ is free in the sense that if
$\lambda \in \Lambda$ and $n \in \ZZ^l$ satisfy $\alpha_n(\lambda) =
\lambda$, then $n = \Z{l}$. As in \cite[Section~5]{KP} we may form
the quotient $k$-graph $\Lambda/\ZZ^l$, and \cite[Theorem~5.7]{KP}
shows that $C^*(\Lambda) \times_{\tilde{\alpha}} \ZZ^l$ is stably
isomorphic to $C^*(\Lambda/\alpha)$. In particular, one can deduce
from this, or from direct arguments, that $\Lambda$ is
$\alpha$-cofinal if and only if $\Lambda/\alpha$ is cofinal, and
$\Lambda$ is aperiodic if and only if $\Lambda/\alpha$ is aperiodic.
In particular, when $\alpha$ is free, $\Lambda/\alpha$ is aperiodic
(respectively cofinal) if and only if $\cpg{\Lambda}{\alpha}{l}$ is
aperiodic (respectively cofinal).

If the action $\alpha$ is not free then, as observed on
\cite[page~176]{PQR2}, the natural definition of $\Lambda/\alpha$
need not yield a category:  the obvious candidate for a composition
map is not necessarily well-defined. The approach of
\cite[Theorem~5.7]{KP} therefore cannot be applied to non-free
actions. In particular, we cannot expect, in general, to be able to
study $C^*(\Lambda) \times_{\tilde\alpha} \ZZ^l$ using a quotient
$k$-graph. However, the crossed-product $k$-graph
$\cpg{\Lambda}{\alpha}{l}$ makes sense regardless, and our results
still apply.
\end{rmk}

\begin{proof}[Proof of Theorem~\ref{thm:first simplicity theorem}]
By Theorem~\ref{thm:crossed product}, it suffices to show that
$C^*(\cpg{\Lambda}{\alpha}{l})$ is simple if and only if $\Lambda$ is
both $\alpha$-aperiodic and $\alpha$-cofinal. By Lemmas
\ref{lem:cofinality equivalence}~and~\ref{lem:aperiodicity
equivalence}, it therefore suffices to show that
$C^*(\cpg{\Lambda}{\alpha}{l})$ is simple if and only if
$\cpg{\Lambda}{\alpha}{l}$ is cofinal and has no local periodicity.
Since $\cpg{\Lambda}{\alpha}{l}$ is row-finite and has no sources,
this follows from \cite[Theorem~3.1]{RobSi}.
\end{proof}

\section{$C^*$-algebraic simplicity criteria}\label{sec:C*-simplicity}

In this section, we reinterpret the hypotheses that $\Lambda$ is
$\alpha$-aperiodic, and that $\Lambda$ is $\alpha$-cofinal
$C^*$-algebraically. Specifically, we re-cast these conditions in
terms of the restriction of the induced action $\tilde{\alpha}$ to
the canonical abelian subalgebra $\overline{D} = \clsp\{s_\lambda
s^*_\lambda : \lambda \in \Lambda\}$ of $C^*(\Lambda)$. To do this,
we insist that $\Lambda$ should be \emph{locally finite with no
sources or sinks} in the sense that for each $p \in \NN^k$ and each
$v \in \Lambda^{\Z{k}}$ there the sets $v\Lambda^p$ and $\Lambda^p v$
are both finite and nonempty. The resulting formulation is almost
identical to \cite[Proposition~8.29]{BrownlowePhD}, and many of the
ideas in the proof are drawn from that argument.

\begin{lemma}\label{lem:canonical endos}
Let $\Lambda$ be a locally finite $k$-graph with no sources or sinks.
Let $s_\Lambda = \{s_\lambda : \lambda \in \Lambda\} \subset
C^*(\Lambda)$ denote the universal generating Cuntz-Krieger
$\Lambda$-family. Let $D$ be the $^*$-subalgebra $\lsp\{s_\lambda
s^*_\lambda : \lambda \in \Lambda\}$ of $C^*(\Lambda)$.

For $a \in D$ and $p \in \NN^k$, there are only finitely many
paths $\eta \in \Lambda^p$ such that $s_\eta a s^*_\eta \not=
0$. Define $\Phi_p : D \to D$ by
\begin{equation}\label{eq:Phi_p}
    \Phi_p(a) := \sum_{\eta \in \Lambda^p} s_\eta a s^*_\eta.
\end{equation}
Then $\|\Phi_p(a)\| = \|a\|$ for $a \in D$, and $\Phi_p$ extends to
an endomorphism, also denoted $\Phi_p$ of $\overline{D}$. Moreover,
the map $p \mapsto \Phi_p$ defines an action of $\NN^k$ on
$\overline{D}$ by endomorphisms.
\end{lemma}
\begin{proof}
We first show that $\{\eta \in \Lambda^p : s_\eta a s^*_\eta
\not= 0\}$ is finite. First note that $P_{r(G)} := \sum_{v \in
r(G)} s_v$ is a left-identity for $a$. Relation~(TCK1) ensures
that $s_w P_{r(G)} = 0$ for $w \not\in r(G)$. Relation~(TCK2)
therefore implies that $s_\eta P_{r(G)} = 0$ whenever $\eta \in
\Lambda^p \setminus \Lambda^p r(G)$. In particular, $\{\eta \in
\Lambda^p : s_\eta a s^*_\eta \not= 0\} \subset \Lambda^p
r(G)$. Since $r(G)$ is finite, and since $\Lambda$ is locally
finite, $\Lambda^p r(G)$ itself is finite. Thus $\{\eta \in
\Lambda^p : s_\eta a s^*_\eta \not= 0\}$ is finite as required.

It is clear that $\Phi_p$ is linear and preserves adjoints. To
see that $\Phi_p$ extends to an endomorphism of $\overline{D}$,
it suffices to show that $\|\Phi_p(a)\| = \|a\|$ for each $a
\in D$, and that $\Phi_p(a)\Phi_p(b) = \Phi_p(ab)$ for $a,b \in
D$.

Fix $a \in D$, and write $a = \sum_{\lambda \in F} a_\lambda
s_\lambda s^*_\lambda \in D$ where $F \subset \Lambda$ is
finite. Let $q := \bigvee_{\lambda \in F} d(\lambda)$. As
$\Lambda$ is row-finite with no sources,
\cite[Proposition~B.1]{RSY2} implies that
\begin{equation}\label{eq:oldCKrel}
s_v = \sum_{\lambda \in v\Lambda^q} s_\lambda s^*_\lambda
\quad\text{for all $v\in \Lambda^0$ and $q \in \NN^k$.}
\end{equation}
We may apply~\eqref{eq:oldCKrel}, to each term in $a$ to obtain
a finite set $G \subset \Lambda^q$ and scalars $\{b_\tau : \tau
\in G\}$ such that $a = \sum_{\tau \in G} b_\tau s_\tau
s^*_\tau$. Since the $s_\tau s^*_\tau$ are mutually orthogonal
projections, we have $\|a\| = \max\{|b_\tau| : \tau \in G\}$.
By definition of $\Phi$ and the Cuntz-Krieger relations we have
\[
\Phi_p(a) = \sum_{\tau \in G, \eta \in \Lambda^p r(\tau)}
    b_\tau s_{\eta\tau} s^*_{\eta\tau}.
\]
Since the $s_{\eta \tau} s^*_{\eta\tau}$ are mutually
orthogonal, $\|\Phi_p(a)\| = \max\{|b_\tau| : \tau \in G,
\Lambda^p r(\tau) \not= \emptyset\}$. Since $\Lambda$ has no
sinks, this gives $\|\Phi_p(a)\| = \|a\|$.

To see that $\Phi_p$ is multiplicative, fix $a, b$ in $D$ and
use~\eqref{eq:oldCKrel} to express $a = \sum_{\tau \in F} a_\tau
s_\tau s^*_\tau$ and $b = \sum_{\rho \in G} b_\rho s_\rho s^*_\rho$
where $F,G$ are finite subsets of $\Lambda^q$ for some fixed $q \in
\NN^k$.
\begin{align}
    \Phi_p(a)\Phi_p(b)
        &= \sum_{\substack{\tau \in F, \eta \in \Lambda^p r(\tau)
        \\
            \rho \in G, \zeta \in \Lambda^p r(\rho)}}
            s_\eta (a_\tau s_\tau s^*_\tau) s^*_\eta
            s_\zeta (b_\rho s_\rho s^*_\rho) s^*_\zeta \nonumber \\
        &= \sum_{\substack{\tau \in F, \eta \in \Lambda^p r(\tau)
        \\
            \rho \in G, \zeta \in \Lambda^p r(\rho)}}
            s_\eta (a_\tau s_\tau s^*_{\eta\tau} s_{\zeta\rho}
            b_\rho s^*_\rho) s^*_\zeta. \label{eq:2ndlastline}
\end{align}
Consider a product $s^*_{\eta\tau} s_{\zeta\rho}$ occurring in
a term in~\eqref{eq:2ndlastline}. Since $F,G \subset
\Lambda^q$, \eqref{eq:oldCKrel} ensures that if the term
$s^*_{\eta\tau} s_{\zeta\rho}$ is nonzero, then $\eta\tau =
\zeta\rho$. Since $d(\eta) = d(\zeta) = p$, the factorisation
property guarantees that $\eta\tau = \zeta\rho$ if and only if
$\eta = \zeta$ and $\tau = \rho$. Hence
\begin{equation}\label{eq:lastline}
\Phi_p(a)\Phi_p(b)
= \sum_{\substack{\tau \in F \cap G, \eta \in \Lambda^p(r(\tau))}}
s_{\eta} (a_\tau b_\tau s^*_\tau s_\tau) s^*_{\eta}.
\end{equation}
To see that~\eqref{eq:lastline} is equal to $\Phi_p(ab)$, we
calculate
\[
ab = \sum_{\tau \in F, \rho \in G} a_\tau b_\rho s_\tau s^*_\tau
s_\rho s^*_\rho = \sum_{\tau \in F \cap G} a_\tau b_\tau s_\tau
s^*_\tau
\]
by~\eqref{eq:oldCKrel}. Applying the formula for $\Phi_p$ to this
expression, we obtain the right-hand side of~\eqref{eq:lastline}.

It remains to show that $p \mapsto \Phi_p$ determines an action of
$\NN^k$; that is, we must show that $\Phi_p\circ \Phi_q = \Phi_{p+q}$
for all $p,q \in \NN^k$. By linearity, it suffices to show that
$\Phi_{p}(\Phi_q(s_\lambda s^*_\lambda)) = \Phi_{p+q}(s_\lambda
s^*_\lambda)$ for all $\lambda \in \Lambda$. Fix $\lambda \in
\Lambda$. Then
\begin{align*}
\Phi_p(\Phi_q(s_\lambda s^*_\lambda))
    &= \Phi_p\Big(\sum_{\eta \in \Lambda^q r(\lambda)}
        s_{\eta\lambda} s^*_{\eta\lambda}\Big)
    &= \sum_{\eta \in \Lambda^q r(\lambda), \zeta \in \Lambda^p r(\eta)}
        s_{\zeta\eta\lambda} s^*_{\zeta\eta\lambda}.
\end{align*}
The factorisation property implies that $(\zeta,\eta) \mapsto
\zeta\eta$ is a bijection of $\{(\zeta,\eta) : \eta \in \Lambda^q
r(\lambda), \zeta \in \Lambda^p r(\eta)\}$ onto $\Lambda^{q+p}
r(\lambda)$. Hence
\[
\Phi_p(\Phi_q(s_\lambda s^*_\lambda))
= \sum_{\xi \in \Lambda^{q+p} r(\lambda)} s_{\xi\lambda} s^*_{\xi\lambda}
= \Phi_{q+p}(s_\lambda s^*_\lambda),
\]
and $\Phi$ is an action as claimed.
\end{proof}

\begin{prop}\label{prp:induced homeos}
Let $\Lambda$ be a locally finite $k$-graph with no sources or sinks,
and let $\overline{D} = \clsp\{s_\lambda s^*_\lambda : \lambda \in
\Lambda\}$. There is a unique isomorphism $\psi : \overline{D} \to
C_0(\Lambda^\infty)$ which takes $s_\lambda s^*_\lambda$ to the
indicator function $\indicator_{\lambda\Lambda^\infty}$. For $p \in
\NN^k$ and $f \in C_0(\Lambda^\infty)$, the endomorphism $\Phi_p$ of
$\overline{D}$ obtained from Lemma~\ref{lem:canonical endos}
satisfies $\psi \circ \Phi_p \circ \psi^{-1}(f) = f\circ \sigma^p$ as
elements of $C_0(\Lambda^\infty)$. If $\alpha$ is an action of
$\ZZ^l$ by automorphisms of $\Lambda$, then for $m \in \ZZ^l$ and $f
\in C_0(\Lambda^\infty)$, the automorphism $\tilde{\alpha}_m$ of
$\overline{D}$ satisfies $\psi \circ \tilde{\alpha}_m \circ
\psi^{-1}(f) = f\circ \alpha^\infty_{-m}$.
\end{prop}
\begin{proof}
The existence of the isomorphism $\psi$ follows from
\cite[Corollary~3.5(i)]{KP}; $\psi$ is unique because
$\overline{D}$ is generated as a $C^*$-algebra by the
$s_\lambda s^*_\lambda$. For the last assertion it suffices to
show that for $\lambda \in \Lambda$
\[
\psi \circ \Phi_p \circ
\psi^{-1}(\indicator_{\lambda\Lambda^\infty})(x) =
\indicator_{\lambda\Lambda^\infty} \circ \sigma^p(x)\quad\text{for
$x
\in
\Lambda^\infty$},
\]
and similarly for $\tilde{\alpha}_m$ and $\alpha^\infty_m$. Fix $x
\in \Lambda^\infty$. Since
\begin{align*}
\psi \circ \Phi_p \circ
\psi^{-1}(\indicator_{\lambda\Lambda^\infty})
    &=\psi(\Phi_p(s_\lambda s^*_\lambda))\\
    &=\psi\Big(\sum_{\xi\in r(\lambda)\Lambda^p}s_{\xi\lambda}
    s^*_{\xi\lambda}\Big)\\
    &=\sum_{\xi\in r(\lambda)\Lambda^p}\indicator_{\xi\lambda},
\end{align*}
we may calculate
\begin{align*}
\psi \circ \Phi_p \circ
\psi^{-1}(\indicator_{\lambda\Lambda^\infty})(x)
    &=\Big(\sum_{\xi\in
    r(\lambda)\Lambda^p}\indicator_{\xi\lambda}\Big)(x)\\
    &=\begin{cases}
        1& \text{if } x(\Z{k}, p+d(\lambda))=\xi\lambda \text{ for
        some } \xi\in r(\lambda)\Lambda^p\\
        0& \text{otherwise}
    \end{cases}\\
    &=\begin{cases}
        1& \text{if } \sigma^p (x)(\Z{k}, d(\lambda))=\lambda \\
        0& \text{otherwise}
        \end{cases}\\
    &=\indicator_{\lambda\Lambda^\infty}(\sigma^p(x)).
\end{align*}
A similar argument establishes the identity involving
$\tilde{\alpha}_m$ and $\alpha^\infty_{-m}$.
\end{proof}

Suppose that $\Lambda$ is a locally finite $k$-graph with no sources
or sinks, and that $\alpha$ is an action of $\ZZ^l$ on $\Lambda$ by
automorphisms. For $(p,m) \in \NN^k \times \NN^l$, the map
$\sigma^p\circ\alpha^\infty_{-m}$ is a local homeomorphism of
$\Lambda^\infty$. We denote this local homeomorphism by
$\tau_{p,m}^{\sigma,\alpha}$. Then $\tau^{\sigma,\alpha} : (p,m)
\mapsto \tau_{p,m}^{\sigma,\alpha}$ is an action of the semigroup
$\NN^k \times \NN^l$ by local homeomorphisms of $\Lambda^\infty$.

Let $X$ be a topological space, and let $\tau$ be an action of
a semigroup $S$ by local homeomorphisms of $X$. As in
\cite{ArchSp, BrownlowePhD},
\begin{Enumerate}
\item We say that the system $(X,\tau)$ is \emph{topologically free}
    if for every pair of distinct elements $s,t \in S$, the set $\{x
    \in X : \tau_s(x) = \tau_t(x)\}$ has empty interior.
\item We say that $x,y \in X$ are \emph{trajectory equivalent} if
    there exist $s,t \in S$ such that $\tau_s(x) = \tau_t(y)$.
\item We say that  $W \subset X$ is \emph{invariant} if $y \in W$ and
    $x$ trajectory equivalent to $y$ imply $x \in W$.
\item We say that $\tau$ is \emph{irreducible} if the only open
    invariant subsets of $X$ are $\emptyset$ and $X$.
\end{Enumerate}

\begin{rmk}\label{rmk:trajectory equivalence}
Suppose that $S = \NN^k \times \NN^l$. Then trajectory equivalence is
an equivalence relation: it is clearly reflexive and symmetric, and
to see that it is transitive, suppose that $x,y,z \in X$, $p,q,p',q'
\in \NN^k$ and $m,n,m',n' \in \NN^l$ satisfy $\tau_{(p,m)}(x) =
\tau_{(q,n)}(y)$ and $\tau_{(p',m')}(y) = \tau_{(q',n')}(z)$. Then
\begin{align*}
\tau_{(p + (q \vee q') - q, m + (n \vee n') - n)}(x)
    &= \tau_{(q \vee q') - q, (n \vee n') - n}(\tau_{(p,m)}(x)) \\
    &= \tau_{(q \vee q') - q, (n \vee n') - n}(\tau_{(q,n)}(y)) \\
    &= \tau_{(q,m) \vee (q',n')}(y).
\end{align*}
Symmetrically, $\tau_{(p + (q \vee q') - q', m + (n \vee n') -
n')}(z) = \tau_{(q,m) \vee (q',n')}(y)$, so that
\[
\tau_{(p + (q \vee q') - q, m + (n \vee n') - n)}(x) = \tau_{(p + (q
\vee q') - q', m + (n \vee n') - n')}(z).
\]
In particular, a set $U \subset X$ is invariant if and only if its
complement $X \setminus U \subset X$ is invariant, so $\tau$ is
irreducible if and only if the only closed invariant subsets of $X$
are $\emptyset$ and $X$.
\end{rmk}

\begin{theorem}\label{thm:simplicity theorem}
Let $\Lambda$ be a locally finite $k$-graph with no sources or sinks,
and let $\alpha$ be an action of $\ZZ^l$ by automorphisms of
$\Lambda$. Then
\begin{enumerate}
\item Every nontrivial ideal $I$ of $C^*(\Lambda)
    \times_{\tilde{\alpha}} \ZZ^l$ satisfies $I \cap
    \pi_{C^*(\Lambda)}(\overline{D}) \not= \{0\}$ if and only if
    $(\Lambda^\infty, \tau^{\sigma, \alpha})$ is topologically
    free.
\item The
ideals $I(a)$ in $C^*(\Lambda) \times_{\tilde{\alpha}} \ZZ^l$
generated by nonzero elements $a$ of
$\pi_{C^*(\Lambda)}(\overline{D})$ are all equal to $C^*(\Lambda)
\times_{\tilde{\alpha}} \ZZ^l$ if and only if
$\tau^{\sigma,\alpha}$ is irreducible.
\end{enumerate}
In particular $C^*(\Lambda) \times_{\tilde{\alpha}} \ZZ^l$ is simple
if and only if $(\Lambda^\infty, \tau^{\sigma,\alpha})$ is
topologically free and $\tau^{\sigma, \alpha}$ is irreducible.
\end{theorem}

\begin{rmk}
The above theorem applies when $l = 0$ so that $\alpha$ is the
trivial action of the trivial group $\{0\}$. In this case,
$\tau^{\sigma,\alpha}$ is just the action $\sigma$ of $\NN^k$
on $\Lambda^\infty$ by shift maps, and we obtain a parallel
result to \cite[Proposition~8.29]{BrownlowePhD} for
locally-finite $k$-graphs.
\end{rmk}

To prove Theorem~\ref{thm:simplicity theorem} we establish two
lemmas. The first establishes that topological freeness is equivalent
to $\alpha$-aperiodicity, and the second that irreducibility is
equivalent to $\alpha$-cofinality. We then apply
Theorem~\ref{thm:first simplicity theorem} to obtain the result.

\begin{lemma}\label{lem:alpha-aper/top free}
Let $\Lambda$ be a locally finite $k$-graph with no sources or sinks,
and let $\alpha$ be an action of $\ZZ^l$ by automorphisms of
$\Lambda$. Then $\Lambda$ is $\alpha$-aperiodic if and only if
$(\Lambda^\infty, \tau^{\sigma, \alpha})$ is topologically free.
\end{lemma}
\begin{proof}
First suppose that $(\Lambda^\infty, \tau^{\sigma, \alpha})$ is
topologically free. Fix $v \in \Lambda^{\Z{k}}$ and $(p,m) \not=
(q,n) \in \NN^k \times \NN^l$. Since $v\Lambda^\infty$ is open in
$\Lambda^\infty$, topological freeness ensures that there exists an
$x \in v\Lambda^\infty$ such that $\tau^{\sigma, \alpha}_{p,m}(x)
\not= \tau^{\sigma, \alpha}_{q,n}(x)$; that is
$\sigma^p(\alpha^\infty_{-m}(x)) \not=
\sigma^q(\alpha^\infty_{-n}(x))$. Since $v, (p,m)$ and $(q,n)$ were
arbitrary, it follows that $\Lambda$ is $\alpha$-aperiodic.

Now suppose that $\Lambda$ is $\alpha$-aperiodic. We must show that
$(\Lambda^\infty, \tau^{\sigma, \alpha})$ is topologically free. Fix
$(p,m) \not= (q,n) \in \NN^k \times \NN^l$. We must show that $\{x
\in \Lambda^\infty : \sigma^p(\alpha^\infty_{-m}(x)) =
\sigma^q(\alpha^\infty_{-n}(x))\}$ has empty interior. Since the sets
$\lambda \Lambda^\infty, \lambda \in \Lambda$ form a basis for the
topology on $\Lambda^\infty$, it suffices to show that for each fixed
$\lambda \Lambda^\infty$, there exists $x \in \lambda\Lambda^\infty$
such that $\sigma^p(\alpha^\infty_{-m}(x)) \not=
\sigma^q(\alpha^\infty_{-n}(x))$. Fix $\lambda \in \Lambda$. By
$\alpha$-aperiodicity, there exists $y \in s(\lambda)\Lambda^\infty$
such that $\sigma^p(\alpha^\infty_{-m}(y)) \not=
\sigma^q(\alpha^\infty_{-n}(y))$, and then $x := \lambda y$ has the
desired property.
\end{proof}

\begin{lemma}\label{lem:alpha-cofinal/irred}
Let $\Lambda$ be a locally finite $k$-graph with no sources or sinks,
and let $\alpha$ be an action of $\ZZ^l$ by automorphisms of
$\Lambda$. Then $\Lambda$ is $\alpha$-cofinal if and only if
$\tau^{\sigma, \alpha}$ is irreducible.
\end{lemma}
\begin{proof}
We follow the proof of \cite[Lemma~8.31]{BrownlowePhD} quite
closely. First suppose that $\Lambda$ is $\alpha$-cofinal. Let
$U$ be a nonempty open invariant subset of $\Lambda^\infty$; we
must show that $U = \Lambda^\infty$. Since $U$ is open and
nonempty, there exists $\lambda \in \Lambda$ such that
$\lambda\Lambda^\infty \subset U$. Fix $x \in \Lambda^\infty$.
Since $\Lambda$ is $\alpha$-cofinal, there exist $p \in \NN^k$,
$m,n \in \NN^l$ and $\mu \in \Lambda$ such that $r(\mu) =
\alpha_{-m}(s(\lambda))$ and $s(\mu) = \alpha_{-n}(x(p))$. We
have $y := \lambda \alpha_m(\mu)
\alpha^\infty_{m-n}(\sigma^p(x)) \in \lambda\Lambda^\infty
\subset U$. Moreover,
\[
\tau^{\sigma,\alpha}_{d(\lambda) + d(\mu),\, n}(y) =
\alpha^\infty_n(\alpha^\infty_{m-n}(\sigma^p(x)) =
\tau^{\sigma,\alpha}_{p,m}(x),
\]
so that $x$ and $y$ are trajectory equivalent. Since $U$ is
invariant, this forces $x \in U$. Since $x \in \Lambda^\infty$ was
arbitrary, it follows that $U = \Lambda^\infty$.

Now suppose that $\Lambda$ is not $\alpha$-cofinal, and fix $v \in
\Lambda^{\Z{k}}$ and $x \in \Lambda^\infty$ such that $\alpha_{-m}(v)
\Lambda \alpha_{-n}(x(p)) = \emptyset$ for all $p \in \NN^k$ and $m,n
\in \NN^l$. We will show that $\tau^{\sigma,\alpha}$ is not
irreducible by constructing an open invariant set $U$ which is equal
to neither $\Lambda^\infty$ nor $\emptyset$. Let
\[
U := \{y \in \Lambda^\infty : \alpha_{-m}(v) \Lambda
\alpha_{-n}(y(p)) \not= \emptyset\text{ for some } p \in \NN^k\text{
and }m,n \in \NN^l\}.
\]
Since $v\Lambda^\infty \not= \emptyset$, we have $U \not= \emptyset$,
and since $x \not\in U$ by construction, we have $U \not=
\Lambda^\infty$.

We claim that $U$ is open. To see this, fix $z \in U$, and let $p \in
\NN^k$ and $m,n \in \NN^l$ satisfy
$\alpha_{-m}(v)\Lambda\alpha_{-n}(z(p)) \not= \emptyset$. Let
$\lambda := z(\Z{k},p)$. For any $z' \in \lambda\Lambda^\infty$ we
have $\alpha_{-n}(z'(p)) = \alpha_{-n}(s(\lambda)) =
\alpha_{-n}(z(p))$, and hence $\lambda\Lambda^\infty \subset U$.
Since $\lambda\Lambda^\infty$ is an open neighbourhood of $z$, it
follows that $U$ is open.

Finally, we claim that $U$ is invariant. Suppose that $y \in U$ and
that $z$ is trajectory equivalent to $y$. Fix $p \in \NN^k$ and $m,n
\in \NN^l$ such that $\alpha_{-m}(v) \Lambda \alpha_{-n}(y(p)) \not=
\emptyset$. Then
\begin{equation}\label{eq:extend}
\alpha_{-(m+h)}(v) \Lambda \alpha_{-(n+h)}(y(p+l)) \not=
\emptyset\text{ for all $l \in \NN^k, h \in \NN^l$}.
\end{equation}
Since $y$ is trajectory equivalent to $z$, there exist $c,d \in
\NN^k$ and $a,b \in \NN^l$ such that $\tau^{\sigma,\alpha}_{c,a}(y) =
\tau^{\sigma,\alpha}_{d,b}(z)$; that is, $\alpha_{-a}(y(c)) =
\alpha_{-b}(z(d))$, and we may assume without loss of generality that
$a \ge n$ and $c \ge p$; say $a = n + h$ and $c = p + l$.
By~\eqref{eq:extend}, we then have
\[
\emptyset \not= \alpha_{-(m + h)}(v) \Lambda \alpha_{-(n +
h)}(y(p+l)) = \alpha_{-(m + h)}(v) \Lambda \alpha_{-a}(y(c)),
\]
and as $\alpha_{-a}(y(c)) = \alpha_{-b}(z(d))$ by choice of
$a,b,c,d$, it follows that $z \in U$. Hence $U$ is invariant, and the
proof is complete.
\end{proof}

\begin{proof}[Proof of Theorem~\ref{thm:simplicity theorem}]
The last assertion follows from (1)~and~(2). By definition of
$\overline{D}$, every ideal of $C^*(\Lambda)$ which intersects
$\overline{D}$ must contain a vertex projection. Hence it suffices to
show: (a) that every ideal of $C^*(\Lambda) \times_{\tilde{\alpha}}
\ZZ^l$ contains a vertex projection if and only if $(\Lambda^\infty,
\tau^{\sigma, \alpha})$ is topologically free; and (b) that the
ideals of $C^*(\Lambda) \times_{\tilde{\alpha}} \ZZ^l$ generated by
the projections $\pi_{C^*(\Lambda)}(s_v)$ are all equal to the whole
crossed product if and only if $\tau^{\sigma,\alpha}$ is irreducible.
But~(a) follows from Lemmas \ref{lem:alpha-aper/top
free}~and~\ref{lem:aperiodicity equivalence} together with
\cite[Proposition~3.6]{RobSi} and~(b) follows from Lemmas
\ref{lem:alpha-cofinal/irred}~and~\ref{lem:cofinality equivalence}
together with \cite[Proposition~3.5]{RobSi}.
\end{proof}

\section{$K$-theory}\label{sec:K-theory}
In this section we consider an action of $\ZZ$ on a row-finite
$1$-graph $E$ with no sources such that either $K_0(C^*(E)) = \{0\}$
or $K_1(C^*(E)) = \{0\}$. In this case, $C^*(\cpg{E}{\alpha}{})$ is a
crossed product of $C^*(E)$ by $\ZZ$ and we can use the
Pimsner-Voiculescu exact sequence to investigate its $K$-theory.

Our main application is to the $2$-graphs discussed in
\cite[Section~3]{PRRS}. These can be realised, using
Proposition~\ref{prp:when cpgraph}, as $2$-graphs of the form
$\cpg{E}{\alpha}{}$ where $E$ is a $1$-graph with no cycles. This
guarantees that $C^*(E)$ is an AF algebra, so has trivial
$K_1$-group, and the Pimsner-Voiculescu sequence provides a
relatively straightforward calculation of the $K$-theory of the
crossed-product algebra $C^*(E) \times_{\tilde{\alpha}} \ZZ$. This
approach is significantly more efficient than the calculations of
\cite[Section~4]{PRRS}, even for the smaller class of rank-2 Bratteli
diagrams considered there.

To state the main result, we need some notation and definitions.

Given a set $X$, we let $\ZZ X$ denote the collection of finitely
supported functions $f : X \to \ZZ$, and regard it as a group under
pointwise addition. We write $\{\delta_x : x \in X\}$ for the
canonical basis for $\ZZ X$.

Let $E$ be a row-finite $1$-graph with no sources. Let $M_E$
denote the connectivity matrix of $E$ given by $M_E(v,w) = |\{e
\in E^1 : r(e) = v, s(e) = w\}|$. We regard $M_E$ as a
homomorphism of $\ZZ E^0$ (implemented by matrix
multiplication).

Given an automorphism $\alpha$ of $E$, we write $\alpha_*$ for the
induced homomorphism $\alpha_* : \ZZ E^0 \to \ZZ E^0$ determined by
$\alpha_*(f)(v) = f(\alpha(v))$. Equivalently, $\alpha_*(\delta_v) =
\delta_{\alpha^{-1}(v)}$ for all $v \in E^0$.

\begin{theorem}\label{thm:induced K-maps for alpha}
Let $E$ be a row-finite $1$-graph with no sources, and let $\alpha$
be an automorphism of $E$. Resume the notation outlined above. Then
$\alpha_*$ commutes with $M_E^t$ and induces an automorphism
$\widetilde{\alpha_*}$ of $\coker(1 - M_E^t)$ satisfying
\[
\widetilde{\alpha_*}(f + \im(1 - M_E^t)) := \alpha_*(f) + \im(1 -
M_E^t).
\]
Furthermore, $\alpha_*$ restricts to an automorphism $\alpha_*|$ of
$\ker(1 - M_E^t)$.

There is an isomorphism $\phi_0 : K_0(C^*(E)) \to \coker(1-M_E^t)$
which satisfies $\phi_0([s_v]) = \delta_v + \im(1 - M_E^t)$ and there
is an isomorphism $\phi_1 : K_1(C^*(E)) \to \ker(1 - M_E^t)$ such
that the diagrams
\[
\begin{CD}
K_0(C^*(E)) @>{\phi_0}>> \coker(1 - M_E^t) \\
@VV{K_0(\tilde\alpha)}V @V{\widetilde{\alpha_*}}VV \\
K_0(C^*(E)) @>{\phi_0}>> \coker(1 - M_E^t)
\end{CD}
\qquad\text{and}\qquad
\begin{CD}
K_1(C^*(E)) @>{\phi_1}>> \ker(1 - M_E^t) \\
@VV{K_1(\tilde\alpha)}V @V{\alpha_*|}VV \\
K_1(C^*(E)) @>{\phi_1}>> \ker(1 - M_E^t)
\end{CD}
\]
commute.
\end{theorem}
\begin{proof}
To see that $\alpha_*$ commutes with $M_E^t$, fix a generator
$\delta_v$ of $\ZZ E^0$ calculate:
\[\begin{split}
\alpha_*(M_E^t \delta_v)
 = \alpha_*\Big(\sum_{r(e) = v} \delta_{s(e)}\Big)
 &= \sum_{r(e) = v} \delta_{\alpha^{-1}(s(e))} \\
 &= \sum_{r(e') = \alpha^{-1}(v)} \delta_{s(e')}
 = M_E^t\delta_{\alpha^{-1}(v)}
 = M_E^t(\alpha_*(\delta_v)).
\end{split}\]

The remaining statements follow from the definitions of the maps
$\tilde{\alpha}$ and $\widetilde{\alpha_*}$, the $K$-theory
calculations for graph algebras of \cite{PR, RSz}, and the naturality
of the Pimsner-Voiculescu exact sequence (see for example
\cite[Section~3]{RSz}).
\end{proof}

\begin{cor}\label{cor:cpg K-th}
Resume the notation of Theorem~\ref{thm:induced K-maps for alpha}.
\begin{Enumerate}
\item
Suppose that $K_1(C^*(E)) = \{0\}$. Then $K_0(C^*(\cpg{E}{\alpha}{}))
\cong \coker(1 - \widetilde{\alpha_*})$ and
$K_1(C^*(\cpg{E}{\alpha}{})) \cong \ker(1 - \widetilde{\alpha_*})$.
\item
Suppose that $K_0(C^*(E)) = \{0\}$. Then $K_0(C^*(\cpg{E}{\alpha}{}))
\cong \ker(1 - \alpha_*|_{\ker(1 - M_E^t)})$ and
$K_1(C^*(\cpg{E}{\alpha}{})) \cong \coker(1 - \alpha_*|_{\ker(1 -
M_E^t)})$.
\end{Enumerate}
\end{cor}
\begin{proof}
The result follows immediately from the Pimsner-Voiculescu exact
sequence for the action $\tilde{\alpha}$ (see
\cite[Theorem~2.4]{PV}).
\end{proof}

\begin{notation}\label{ntn:A,B}
Let $E$ be a row-finite $1$-graph with no sources, and let $\alpha$
be an action of $\ZZ$ on $E$ by automorphisms such that the orbit of
each vertex under $\alpha$ is finite. For each $v \in E^0$, let $C(v)
:= \{\alpha_n(v) : n \in \ZZ\}$ be the orbit of $v$. Let $\mathcal{C}
:= \{C(v) : v \in E^0\}$ be the collection of all orbits of vertices
under $\alpha$. Define integer-valued $\mathcal{C} \times
\mathcal{C}$ matrices $A$ and $B$ as follows: for $C_1, C_2 \in
\mathcal{C}$, define
\[
A_{C_1,C_2} := |C_1 E^1 C_2|/|C_1|
 \quad\text{and}\quad
B_{C_1,C_2} := |C_1 E^1 C_2|/|C_2|.
\]
An argument like the proof of \cite[Lemma~4.2]{PRRS}, shows that for
any $v_1 \in C_1$ and $v_2 \in C_2$, we have $A_{C_1, C_2} = |v_1 E^1
C_2|$ and $B_{C_1, C_2} = |C_1 E^1 v_2|$. We regard $A$ and $B$ as
homomorphisms of $\ZZ\mathcal{C}$ regarded as a group under addition.
\end{notation}

\begin{prop}\label{prp:genPRRS}
Let $E$ be a row-finite $1$-graph with no sources, and let $\alpha$
be an action of $\ZZ$ on $E$ by automorphisms. Suppose that the orbit
of each vertex $v \in E^0$ under $\alpha$ is finite, and let
$\mathcal{C}$, and $A,B : \ZZ \mathcal{C} \to \ZZ \mathcal{C}$ be as
in Notation~\ref{ntn:A,B}.
\begin{Enumerate}
\item If $K_1(C^*(E)) = \{0\}$, then
\[
K_0(C^*(\cpg{E}{\alpha}{})) \cong \coker(1 - A^t)
 \quad\text{and}\quad
K_1(C^*(\cpg{E}{\alpha}{})) \cong \coker(1 - B^t).
\]
\item If $K_0(C^*(E)) = \{0\}$, then
\[
K_0(C^*(\cpg{E}{\alpha}{})) \cong \ker(1 - B^t)
 \quad\text{and}\quad
K_1(C^*(\cpg{E}{\alpha}{})) \cong \ker(1 - A^t).
\]
\end{Enumerate}
\end{prop}
\begin{proof}
Let $M_E$ be the adjacency matrix of $E$. Let $\phi := 1 - M_E^t :
\ZZ E^0 \to \ZZ E^0$, and let $\psi := 1 - \alpha_* : \ZZ E^0 \to \ZZ
E^0$.

The strategy is to define maps $\phi|$, $\tilde{\phi}$,
$\tilde{\psi}$, $q_{\phi}$, $q_{\psi}$, $q_{\phi|}$,
$q_{\tilde{\phi}}$, $q_{\psi}|$ and $\widetilde{q_{\psi}}$ which make
the 16-term diagram illustrated in Figure~\ref{fig:16-term} commute.
\begin{figure}[ht]
\[
\begin{tikzpicture}
\newlength{\myinnersep}
\setlength{\myinnersep}{4pt}
% \node[inner sep=\myinnersep] (tttlll) at (0,7.5) {0};
 \node[inner sep=\myinnersep] (tttll) at (3,7.5) {0};
 \node[inner sep=\myinnersep] (tttl) at (6,7.5) {0};
 \node[inner sep=\myinnersep] (tttr) at (9,7.5) {0};
 \node[inner sep=\myinnersep] (tttrr) at (12,7.5) {0};
% \node[inner sep=\myinnersep] (tttrrr) at (15,7.5) {0};
%
 \node[inner sep=\myinnersep] (ttlll) at (0.5,6) {0};
 \node[inner sep=\myinnersep] (ttll) at (3,6) {$\ker(\phi) \cap
 \ker(\psi)$};
 \node[inner sep=\myinnersep] (ttl) at (6,6) {$\ker(\phi)$};
 \node[inner sep=\myinnersep] (ttr) at (9,6) {$\ker(\phi)$};
 \node[inner sep=\myinnersep] (ttrr) at (12,6)
 {$\ker(\tilde{\phi})$};
 \node[inner sep=\myinnersep] (ttrrr) at (14.5,6) {0};
 \node[inner sep=\myinnersep] (tlll) at (0.5,4.5) {0};
 \node[inner sep=\myinnersep] (tll) at (3,4.5) {$\ker(\psi)$};
 \node[inner sep=\myinnersep] (tl) at (6,4.5) {$\ZZ E^0$};
 \node[inner sep=\myinnersep] (tr) at (9,4.5) {$\ZZ E^0$};
 \node[inner sep=\myinnersep] (trr) at (12,4.5) {$\coker(\psi)$};
 \node[inner sep=\myinnersep] (trrr) at (14.5,4.5) {0};
 \node[inner sep=\myinnersep] (blll) at (0.5,3) {0};
 \node[inner sep=\myinnersep] (bll) at (3,3) {$\ker(\psi)$};
 \node[inner sep=\myinnersep] (bl) at (6,3) {$\ZZ E^0$};
 \node[inner sep=\myinnersep] (br) at (9,3) {$\ZZ E^0$};
 \node[inner sep=\myinnersep] (brr) at (12,3) {$\coker(\psi)$};
 \node[inner sep=\myinnersep] (brrr) at (14.5,3) {0};
 \node[inner sep=\myinnersep] (bblll) at (0.5,1.5) {0};
 \node[inner sep=\myinnersep] (bbll) at (3,1.5) {$\coker(\phi|)$};
 \node[inner sep=\myinnersep] (bbl) at (6,1.5) {$\coker(\phi)$};
 \node[inner sep=\myinnersep] (bbr) at (9,1.5) {$\coker(\phi)$};
 \node[inner sep=\myinnersep] (bbrr) at (12,1.5)
 {$\coker(\tilde{\phi})$};
 \node[inner sep=\myinnersep] (bbrrr) at (14.5,1.5) {0};
 \node[inner sep=\myinnersep] (bbbll) at (3,0) {0};
 \node[inner sep=\myinnersep] (bbbl) at (6,0) {0};
 \node[inner sep=\myinnersep] (bbbr) at (9,0) {0};
 \node[inner sep=\myinnersep] (bbbrr) at (12,0) {0};
 \draw[-latex] (tttll.south) -- (ttll.north);
 \draw[-latex] (tttl.south) -- (ttl.north);
 \draw[-latex] (tttr.south) -- (ttr.north);
 \draw[-latex] (tttrr.south) -- (ttrr.north);
 \draw[right hook-latex] (ttll.south) -- (tll.north);
 \draw[right hook-latex] (ttl.south) -- (tl.north);
 \draw[right hook-latex] (ttr.south) -- (tr.north);
 \draw[right hook-latex] (ttrr.south) -- (trr.north);
 \draw[-latex] (tll.south) -- node[right]{${}_{\phi|}$} (bll.north);
 \draw[-latex] (tl.south) -- node[right]{${}_\phi$} (bl.north);
 \draw[-latex] (tr.south) -- node[right]{${}_\phi$} (br.north);
 \draw[-latex] (trr.south) -- node[right]{${}_{\tilde{\phi}}$}
 (brr.north);
 \draw[-latex] (bll.south) -- node[right]{${}_{q_{\phi|}}$}
 (bbll.north);
 \draw[-latex] (bl.south) -- node[right]{${}_{q_{\phi}}$}
 (bbl.north);
 \draw[-latex] (br.south) -- node[right]{${}_{q_{\phi}}$}
 (bbr.north);
 \draw[-latex] (brr.south) -- node[right]{${}_{q_{\tilde{\phi}}}$}
 (bbrr.north);
 \draw[-latex] (bbll.south) -- (bbbll.north);
 \draw[-latex] (bbl.south) -- (bbbl.north);
 \draw[-latex] (bbr.south) -- (bbbr.north);
 \draw[-latex] (bbrr.south) -- (bbbrr.north);
 \draw[-latex] (ttlll.east) -- (ttll.west);
 \draw[-latex] (tlll.east) -- (tll.west);
 \draw[-latex] (blll.east) -- (bll.west);
 \draw[-latex] (bblll.east) -- (bbll.west);
 \draw[right hook-latex] (ttll.east) -- (ttl.west);
 \draw[right hook-latex] (tll.east) -- (tl.west);
 \draw[right hook-latex] (bll.east) -- (bl.west);
 \draw[-latex] (bbll.east) -- node[above]{${}_\pi$} (bbl.west);
 \draw[-latex] (ttl.east) -- node[above]{${}_{\psi|}$} (ttr.west);
 \draw[-latex] (tl.east) -- node[above]{${}_\psi$} (tr.west);
 \draw[-latex] (bl.east) -- node[above]{${}_\psi$} (br.west);
 \draw[-latex] (bbl.east) -- node[above]{${}_{\tilde{\psi}}$}
 (bbr.west);
 \draw[-latex] (ttr.east) -- node[above]{${}_{q_\psi|}$}
 (ttrr.west);
 \draw[-latex] (tr.east) -- node[above]{${}_{q_{\psi}}$} (trr.west);
 \draw[-latex] (br.east) -- node[above]{${}_{q_{\psi}}$} (brr.west);
 \draw[-latex] (bbr.east) --
 node[above]{${}_{\widetilde{q_{\psi}}}$}
 (bbrr.west);
 \draw[-latex] (ttrr.east) -- (ttrrr.west);
 \draw[-latex] (trr.east) -- (trrr.west);
 \draw[-latex] (brr.east) -- (brrr.west);
 \draw[-latex] (bbrr.east) -- (bbrrr.west);
\end{tikzpicture}\]
\caption{The 16-term diagram for $\phi$ and $\psi$.}
\label{fig:16-term}
\end{figure}
Under the hypothesis that one of $K_0(C^*(E))$, $K_1(C^*(E))$ is
trivial, we will use the Sixteen Lemma to deduce that all rows and
columns of Figure~\ref{fig:16-term} are exact and hence that
\begin{gather}
 \ker(1 - \widetilde{\alpha_*}) \cong \coker(\phi|),
 \quad\text{and}\quad
 \quad\coker(1 - \widetilde{\alpha_*}) \cong \coker(\tilde{\phi}),
 \label{eq:16 Lemma 1} \\
 \ker((1 - \alpha_*)|_{\ker(\phi)}) \cong \ker(\phi|),
 \quad\text{and}\quad
 \coker((1 - \alpha_*)|_{\ker(\phi)}) \cong \ker(\tilde{\phi}).
 \label{eq:16 Lemma 2}
\end{gather}
Finally, we will establish the existence of isomorphisms $\theta_A :
\coker(\psi) \to \ZZ\mathcal{C}$ and $\theta_B : \ker(\psi) \to
\ZZ\mathcal{C}$ satisfying
\[
\theta_A \circ \tilde\phi = (1 - A^t) \theta_A
 \quad\text{and}\quad
\theta_B \circ \phi| = (1 - B^t) \theta_B.
\]
Combining these with Corollary~\ref{cor:cpg K-th} will complete the
proof.

To define the maps $\phi|$, $\tilde{\phi}$ and $\tilde{\psi}$ in
Figure~\ref{fig:16-term}, recall that Theorem~\ref{thm:induced K-maps
for alpha} implies that $M_E^t$ and $\alpha_*$ commute, and hence
that $\phi$ and $\psi$ commute. It follows that $\phi$ restricts to a
homomorphism $\phi|$ of the kernel of $\psi$ and induces a
homomorphism $\tilde{\phi}$ of $\coker(\psi)$ which satisfies
\[
\tilde{\phi}(f + \psi(\ZZ E^0)) = \phi(f) + \psi(\ZZ E^0).
\]
In a similar fashion, $\psi$ restricts to a homomorphism
$\psi|$ of $\ker(\phi)$ and induces a homomorphism
$\tilde{\psi}$ of $\coker(\phi)$.

The maps $q_{\phi}$, $q_{\psi}$, $q_{\phi|}$ and $q_{\tilde{\phi}}$
in Figure~\ref{fig:16-term} are the natural quotient maps: for
example, $q_\phi$ is defined by $q_\phi(f) = f + \phi(\ZZ E^0)$.

The homomorphism $\pi$ in Figure~\ref{fig:16-term} is defined by
$\pi(f + \phi(\ker(\psi))) := f + \phi(\ZZ E^0)$ for $f \in \ZZ E^0$;
this is well-defined because $\phi(\ker(\psi)) \subset \phi(\ZZ
E^0)$.

The map $q_\psi|$ in Figure~\ref{fig:16-term} is the restriction of
$q_\psi$ to the kernel of $\phi$; this has range in
$\ker(\tilde{\phi})$ by definition of $\tilde{\phi}$.

The map $\widetilde{q_{\psi}}$ in Figure~\ref{fig:16-term} is defined
by $\widetilde{q_\psi}(f + \phi(\ZZ E^0)) := q_{\psi}(f) +
\tilde{\phi}(\coker(\psi))$. To see that this is well-defined, note
that for $f \in \phi(\ZZ E^0)$, we have $f = \phi(g)$ for some $g \in
\ZZ E^0$. Since $q_{\psi}(f) = q_{\psi}(\phi(g)) =
\tilde{\phi}(\psi(g))$ by definition of $q_{\tilde{\phi}}$, it
follows that $q_{\psi}(f) \in\tilde{\phi}(\coker(\psi))$ as required.

We now need to show that the squares in Figure~\ref{fig:16-term}
commute, and that if one of $K_0(C^*(E))$ or $K_1(C^*(E))$ is
trivial, then all rows and columns are exact.

As mentioned above, Theorem~\ref{thm:induced K-maps for alpha} shows
that $\alpha_*$ and $M_E^t$ commute, and it follows that the middle
square of Figure~\ref{fig:16-term} commutes. The other squares
commute by definition of the maps involved.

The middle two rows and all columns are clearly exact. Theorem~4.2.4
of \cite{PR} shows that $\ker(\phi) \cong K_1(C^*(E))$, and
$\coker(\phi) \cong K_0(C^*(E))$. Hence if $K_1(C^*(E)) = \{0\}$,
then the terms in the top row are all equal to $\{0\}$, and that row
is trivially exact, and likewise if $K_0(C^*(E)) = \{0\}$, then the
terms in the bottom row are all equal to $\{0\}$ and that row is
trivially exact.

In either case, we may apply the Sixteen Lemma to deduce that the
remaining row of the diagram is exact. The exactness of the top and
bottom rows establishes the formulae \eqref{eq:16 Lemma
1}~and~\eqref{eq:16 Lemma 2}.

By Corollary~\ref{cor:cpg K-th}, it therefore suffices to show that
$\coker(\tilde{\phi}) \cong \coker(1 - A^t)$, $\coker(\phi|) \cong
\coker(1 - B^t)$, $\ker(\tilde{\phi}) \cong \ker(1 - A^t)$ and
$\ker(\phi|) \cong \ker(1 - B^t)$.

Let $V$ be a subset of $E^0$ which contains exactly one
representative of each orbit $C \in \mathcal{C}$. Then $\coker(\psi)$
is generated by the classes $\{\delta_v + \psi(\ZZ E^0) : v \in V\}$,
and $\delta_v - \delta_w \in \psi(\ZZ E^0)$ if and only if $C(v) =
C(w)$. Hence there is an isomorphism $\theta_A : \coker(\psi) \to
\ZZ\mathcal{C}$ such that $\theta_A(\delta_v + \psi(\ZZ E^0)) =
\delta_{C(v)}$ for all $v \in E^0$.

Now $\tilde{\phi}$ takes $\delta_v + \psi(\ZZ E^0)$ to $\sum_{r(e) =
v} \delta_{s(e)} + \psi(\ZZ E^0)$. Applying $\theta_A$, we have
\begin{align}
\theta_A(\tilde{\phi}(\delta_v + \psi(\ZZ E^0))
 &= \sum_{r(e) = v} \delta_{C(s(e))} \nonumber\\
 &= \sum_{C \in \mathcal{C}} |\{e \in r^{-1}(v) : C =
C(s(e))\}| \delta_C.\label{eq:the C mess}
\end{align}
The expression~\eqref{eq:the C mess} is equal to $(1 -
A^t)(\delta_{C(v)})$ by definition of $A$, and $(1 -
A^t)(\delta_{C(v)}) = (1 - A^t)(\theta_A(\delta_v + \psi(\ZZ E^0)))$
by definition of $\theta_A$. Hence the isomorphism $\theta_A$
intertwines $(1 - A^t)$ and $\tilde{\phi}$, establishing that
$\coker(\tilde{\phi}) \cong \coker(1 - A^t)$ and $\ker(\tilde{\phi})
\cong \ker(1 - A^t)$.

Next note that the map $\alpha_*$ permutes the point-masses
associated to the vertices in each $C \in \mathcal{C}$. Hence
$\ker(\psi)$ is the subgroup of $\ZZ E^0$ generated by
$\{\indicator_C : C \in \mathcal{C}\}$. In particular, there is an
isomorphism $\theta_B : \ker(\psi) \to \ZZ\mathcal{C}$ satisfying
$\indicator_C \mapsto \delta_C$. For $C \in \mathcal{C}$, we have
\begin{equation}\label{eq:coker 1-Mt}
\phi(\indicator_C)
 = \phi(\sum_{v \in C} \delta_v)
 = \sum_{v \in C} \delta_v - \sum_{r(e) = v} \delta_{s(e)}.
\end{equation}
As $v$ ranges over $C$, the vertices $s(e)$ range over all vertices
in orbits $C' \in \mathcal{C}$ such that $r^{-1}(C) \cap s^{-1}(C')
\not= \emptyset$. Moreover, for a fixed vertex $v'$ on in an orbit
$C'$ with $r^{-1}(C) \cap s^{-1}(C') \not= \emptyset$, the term
$-\delta_{v'}$ occurs in the right-hand side of~\eqref{eq:coker 1-Mt}
precisely once for each edge $e$ with $r(e) \in C$ and $s(e) = v'$.
By definition of the matrix $B$, these calculations establish that
$(1 - B^t)\circ\theta_B = \theta_B\circ\phi|_{\ker(\psi)}$. Hence
$\coker(1 - B^t) \cong \coker(\phi|)$ and $\ker(1-B^t) \cong
\ker(\phi|)$ as required.
\end{proof}

\begin{rmk}
Under the hypotheses of Proposition~\ref{prp:genPRRS}, the action
$\alpha$ of $\ZZ$ on $E$ is never free. Consequently, we cannot
calculate the $K$-theory of $C^*(E) \times_{\tilde\alpha} \ZZ$ using
the $K$-theory calculations for graph $C^*$-algebras and the
quotient-graph construction of~\cite{KP1}.
\end{rmk}

\begin{example}
As in Section~3 of \cite{PRRS}, let $\Lambda$ be a row-finite
$2$-graph with no sources such that $\Lambda^{(\NN,0)}$
contains no cycles and each vertex $v \in \Lambda^{\Z{2}}$ is
the range of an isolated cycle $\Lambda^{(0,\NN)}$.

By Proposition~\ref{prp:when cpgraph}, $\Lambda$ is isomorphic to
$\cpg{\Lambda^{(\NN,0)}}{\alpha}{}$ where $\alpha$ is determined by
factorisations through paths in $\Lambda^{(0,\NN)}$. Since
$\Lambda^{(\NN,0)}$ has no cycles, $C^*(\Lambda^{(\NN,0)})$ is an AF
algebra and hence has trivial $K_1$-group. We can therefore apply
Proposition~\ref{prp:genPRRS} to obtain expressions for
$K_*(C^*(\Lambda))$. In particular, our results generalise
\cite[Theorem~4.3(2)]{PRRS} to cover all $2$-graphs described in
Section~3 of \cite{PRRS}.
\end{example}

\end{document}